\documentclass[12pt]{article}
\usepackage[T1]{fontenc}
\usepackage{amsmath,amsthm}
\usepackage[utopia]{mathdesign}
\usepackage[all]{xy}
\usepackage{hyperref}
\usepackage{mathdots}
\usepackage{stmaryrd}

\newcommand{\mr}[1]{\mathrm{#1}}
\newcommand{\mf}[1]{\mathfrak{#1}}
\newcommand{\mc}[1]{\mathcal{#1}}
\newcommand{\mb}[1]{\mathbb{#1}}

\newcommand{\Z}{\mb{Z}}
\newcommand{\Q}{\mb{Q}}
\newcommand{\zp}{\mb{Z}_p}
\newcommand{\qp}{\mb{Q}_p}
\newcommand{\C}{\mb{C}}
\newcommand{\R}{\mb{R}}
\newcommand{\HH}{\mb{H}}

\newcommand{\F}{\mb{F}}
\newcommand{\qbar}{\overline{\Q}}
\newcommand{\Fbar}{\overline{F}}
\newcommand{\qpbar}{\overline{\Q_p\!\!}\,\,}
\newcommand{\qpbt}{\overline{\Q_p\!\!}{\!\!\!\!\!\!\!\phantom{\bar{\Q}}}^{\times}}
\newcommand{\qbt}{\overline{\Q}{\!\!\!\!\!\!\!\phantom{\bar{\Q}}}^{\times}}

\newcommand{\Zp}{\mb{Z}_p}

\newcommand{\La}{\Lambda}
\newcommand{\et}{\text{\'et}}

\newcommand{\ord}{\mr{ord}}
\newcommand{\BM}{\mr{BM}}
\newcommand{\Fr}{\mr{Fr}}

\newcommand{\sub}{\mr{sub}}
\newcommand{\quo}{\mr{quo}}

\newcommand{\Eis}{\mf{m}}
\newcommand{\cusps}{\{\mr{cusps}\}}

\newcommand{\mcS}{\mc{S}}
\newcommand{\tS}{\mc{M}}
\newcommand{\mcT}{\mc{T}}
\newcommand{\tT}{\widetilde{\mc{T}}}
\newcommand{\mcP}{P}

\newcommand{\mbT}{\mb{T}}
\newcommand{\tbT}{\widetilde{\mbT}}
\newcommand{\tGa}{\widetilde{\Gamma}}
\newcommand{\tor}{\mr{tor}}

\newcommand{\tK}{L}
\newcommand{\mcL}{\mc{L}}

\newcommand{\cotimes}{\,\hat{\otimes} \,}
\newcommand{\ps}[1]{\llbracket #1 \rrbracket}
\newcommand{\la}[1]{\lambda_{\frac{#1}{N}}}
\newcommand{\zr}{\zeta_{Np^r}}
\newcommand{\xsim}{\xrightarrow{\smash{\ensuremath{\sim}}}}

\newcommand{\vpr}{\varprojlim_r}

\usepackage{mathtools}
\newlength{\arrow}
\DeclareMathOperator{\Hom}{Hom}
\DeclareMathOperator{\End}{End}
\DeclareMathOperator{\Gal}{Gal}

\DeclareMathOperator{\GL}{GL}
\DeclareMathOperator{\SL}{SL}
\DeclareMathOperator{\PGL}{PGL}
\DeclareMathOperator{\PU}{PU}
\DeclareMathOperator{\SO}{SO}
\DeclareMathOperator{\SU}{SU}
\DeclareMathOperator{\PO}{PO}
\DeclareMathOperator{\reg}{reg}
\DeclareMathOperator{\Spec}{Spec}

\DeclareMathOperator{\Pic}{Pic}

\DeclareMathOperator{\cha}{char}
\DeclareMathOperator{\class}{class}
\DeclareMathOperator{\image}{image}

\DeclareMathOperator{\HS}{HS}

\DeclareMathOperator{\cor}{cor}

\newcommand{\bs}{\operatorname{\backslash}}

\numberwithin{equation}{subsection}

\newtheorem*{thm}{Theorem}
\newtheorem*{conj}{Conjecture}

\theoremstyle{definition}

\newtheorem{sbpara}[equation]{}
\newtheorem*{question}{Question}

\theoremstyle{remark}

\newtheorem*{ack}{Acknowledgments}

\renewcommand{\baselinestretch}{1.2}

\begin{document}

\title{Modular symbols in Iwasawa theory}
\author{Takako Fukaya, Kazuya Kato, and Romyar Sharifi}
\date{}
\maketitle

\section{Introduction}

\begin{sbpara} \label{explicit}

The starting point of this paper is the fascinatingly simple and explicit map
$$
	[u:v] \mapsto \{1-\zeta_N^u,1-\zeta_N^v\}
$$
that relates the worlds of geometry/topology and arithmetic \cite{bus, sh-Lfn}.
Here, for $N \ge 1$,
\begin{enumerate}
	\item[$\bullet$]  $[u:v]$ is a Manin symbol in the relative homology group $H_1(X_1(N),\cusps,\Z)$, 
	and
	\item[$\bullet$] $\{1-\zeta_N^u,1-\zeta_N^v\}$ is a Steinberg symbol in the algebraic $K$-group
	$K_2(\Z[\zeta_N,\tfrac{1}{N}])$, 
\end{enumerate}
where $u, v \in \Z/N\Z$ are nonzero numbers with $(u,v) = (1)$, and
$\zeta_N$ is a primitive $N$th root of unity.

\end{sbpara}

\begin{sbpara} 

The above map connects two different worlds in the following manner:
\begin{eqnarray*}
	\text{geometric theory of } \GL_2 &\Longrightarrow& \text{arithmetic theory of } \GL_1
\end{eqnarray*}
over the field $\Q$.   Here, if we consider the geometry of the modular curve $X_1(N)$ on the left, then we consider the arithmetic of the cyclotomic field $\Q(\zeta_N)$ on the right.  This connection is conjectured to be a correspondence if we work modulo the Eisenstein ideal that is defined in \ref{Eis}:
\begin{eqnarray*}
	\text{geometric theory of } \GL_2 \text{ modulo the Eisenstein ideal} &\Longleftrightarrow& 
	\text{arithmetic theory of } \GL_1.
\end{eqnarray*}
More generally, we are dreaming that there is a strong relationship 
\begin{eqnarray*}
	\text{geometric theory of }\GL_d \text{ modulo the Eisenstein ideal} &\Longleftrightarrow&  
	\text{arithmetic theory of } \GL_{d-1}
\end{eqnarray*}
over global fields.  Our goals are to survey what is known and to explain this dream.

\end{sbpara}

\begin{sbpara} \label{eis}

The connection with the Eisenstein ideal for $\GL_2$ over $\Q$ appears as follows.  The homology group we consider has the action of a Hecke algebra which contains an Eisenstein ideal, and the map of \ref{explicit} factors through the quotient of the homology by this ideal \cite{fk-proof}.  The truth of this is deep and mysterious; it is the idea of specializing at the cusp at $\infty$.  This is the key to the connection between $\GL_2$ and $\GL_1$.  

\end{sbpara}

\begin{sbpara} 

We note that there exist two technical issues with our simple presentation of the ``map'' in \ref{explicit}.  We left out those Manin symbols in which one of $u$ or $v$ is $0$, which are needed to generate the relative homology group.  Also, the map is only well-defined as stated if we first invert $2$ and then project to the fixed part under complex conjugation (see Section \ref{modsym}).

\end{sbpara}

\begin{sbpara}

Let us consider the case that $N$ is a power of an odd prime $p$ and work only with $p$-parts.  Consider the quotient
$$
	P_r = H_1(X_1(p^r),\zp)^+/I_rH_1(X_1(p^r),\zp)^+
$$
of the fixed part of homology under complex conjugation by the action of the Eisenstein ideal $I_r$ in the cuspidal Hecke algebra of weight $2$ and level $p^r$.  By the well-known relationship between $K_2$ and $H^2_{\et}$ of $\Z[\zeta_{p^r},\tfrac{1}{p}]$, the map of \ref{explicit} yields a well-defined map
$$
	\varpi_r \colon P_r \to H^2_{\et}(\Z[\zeta_{p^r},\tfrac{1}{p}],\zp(2))^+
$$
that sends the image of $[u:v]$ in $P_r$ to the cup product $(1-\zeta_{p^r}^u) \cup (1-\zeta_{p^r}^v)$.

\end{sbpara}

\begin{sbpara}

Let us connect this with Iwasawa theory for $\GL_1$.  As we increase $r$, the maps $\varpi_r$ are compatible.   The group $H^2_{\et}(\Z[\zeta_{p^r},\tfrac{1}{p}],\zp(2))^+$ is related to the $p$-part $A_r$ of the class group of $\Q(\zeta_{p^r})$ in the sense that its reduction modulo $p^r$ is isomorphic to the Tate twist of $A_r^-/p^rA_r^-$.   
So, if we let $P = \varprojlim_r P_r$ and $X = \varprojlim_r A_r$, then we obtain a map
$$
	\varpi \colon P \to X^-(1)
$$
that relates geometry of the tower of curves $X_1(p^r)$ modulo the Eisenstein ideal to Iwasawa theory over the
union of cyclotomic fields $\Q(\zeta_{p^r})$.  It is a map of Iwasawa modules under the action of
inverses of diamond operators on the left and of Galois elements on the right.  

\end{sbpara}

\begin{sbpara}

In \cite{sh-Lfn}, the map $\varpi$ is conjectured to be an isomorphism.  If this conjecture is true, we can understand the arithmetic object $X^-$ by using the geometric object $P$.  The Iwasawa main conjecture states that the characteristic ideal of $X^-$ is the equivariant $p$-adic $L$-function.  On the other hand, the characteristic ideal of $P$ under the inverse diamond action can be computed to be a multiple of the Tate twist 
$\xi$ of this $L$-function.  If the characteristic ideals of $X^-(1)$ and $P$ are equal, then the main conjecture follows as a consequence of the analytic class number formula.  Therefore, the conjecture that $\varpi$ is an isomorphism is an explicit refinement of the Iwasawa main conjecture.

\end{sbpara}

\begin{sbpara} In their proof of Iwasawa main conjecture \cite{mw}, Mazur and Wiles, expanding upon the work of Ribet \cite{ribet}, considered the relationship between the geometric theory of $\GL_2$ and the arithmetic theory of $\GL_1$.  Using roughly their methods, we can define a map 
$$
	\Upsilon \colon X^-(1) \to P.
$$
More precisely, $\Upsilon$ is constructed out of the Galois action on the projective limit of the reduction of \'etale homology groups $H_1^{\et}(X_1(p^r)_{/\qbar},\zp)$ modulo the Eisenstein ideal. 
The expectation in \cite{sh-Lfn} is that the maps $\varpi \colon P\to X^-(1)$ and $\Upsilon \colon X^-(1)\to P$ are inverse to each other.  The best evidence we have for this is the equality $\xi' \Upsilon \circ \varpi = \xi'$ after multiplication by the derivative $\xi'$ of $\xi$, which is proven in \cite{fk-proof}.  If the $p$-adic $L$-function $\xi$ has no multiple zeros, this yields the conjecture up to $p$-torsion in $P$.

\end{sbpara}

\begin{sbpara}

An analogous result for the rational function field $\F_q(t)$ over a finite field can be proven by following the analogy between $\F_q(t)$ and $\Q$. In both cases, the key point of the proof is that $(1-\zeta_N^u, 1-\zeta_N^v)$ and its analogue for $\F_q(t)$ are values at the infinity cusp of the ``zeta elements,'' which is to say Beilinson elements and their analogues for $\F_q(t)$, which live in $K_2$ of the modular curve $X_1(N)$ and its Drinfeld analogue for $\F_q(t)$. 

\end{sbpara}

\begin{sbpara}  \label{key_diag}
 
For both $\Q$ and $\F_q(t)$, the philosophy is that $\xi' \Upsilon \circ \varpi$ is the reduction modulo the Eisenstein ideal $I$  of a map involving zeta elements.
Roughly speaking, the proof consists firstly of the demonstration of the existence of a commutative diagram
$$
	\SelectTips{eu}{} \xymatrix{ \mcS \ar[d]_{\mod I} \ar[r]^{z}  & K \ar[d]^{\infty} \ar[r]^{\mr{reg}} & \mf{S} 
	\ar[d]^{\bmod I} \\
	P \ar[r]^-{\varpi} & Y  \ar[r]^-{\xi'\Upsilon} & P. }
$$
Here, $\mcS$ is the space of modular symbols, the  map $z$ takes modular symbols to zeta elements in the $K_2$-group $K$ of a modular curve, $\mf{S}$ is a space of $p$-adic cusp forms, $\mr{reg}$ is the $p$-adic regulator map, and $Y$ is either $X^-(1)$ or its analogue for $\F_q(t)$.  The vertical arrows denoted by ``mod $I$'' are obtained by reduction modulo the Eisenstein ideal $I$ (see Section 2.7 for details), and $\infty$ is given by specialization at a cusp at infinity.  Secondly, it entails a computation of the regulator map on zeta elements that tells us that the composition $\mcS \to K \to \mf{S} \to P$ coincides with $\xi'$ times the projection $\mcS\to P$.
\end{sbpara}

\begin{sbpara}

In this survey paper, we explain the key ideas and concepts of our work, putting aside many of the technical details that must arise in a careful treatment.  While we do our best to strike a balance, the reader should be aware that some of the statements we make require minor modifications in order that they hold.

The structure of the paper is as follows.  In Section 2, we describe the original case of the conjectures and outline the proof of the above result.  In Section 3, we discuss and outline the proof of the analogue for $\F_q(t)$.   In Section 4, we describe what might be expected for $\GL_d$.

\end{sbpara}

\section{The case of $\GL_2$ over $\Q$} \label{GL_2-Q}

Fix an odd prime number $p$ and an integer $N\geq 1$ which is not divisible by $p$.  Let $r \ge 1$, which will vary.  Let $\qbar$ be the algebraic closure of $\Q$ in $\C$, and let us fix an embedding of $\qbar$ in $\qpbar$.

Recall that we want to understand the picture:
\begin{eqnarray*}
	\SelectTips{eu}{} 
	\xymatrix{ \text{geometric theory of } \GL_2 \text{ modulo the Eisenstein ideal} \ar@<0.5ex>[r]^-{\varpi} & 
	\ar@<0.5ex>[l]^-{\Upsilon} 
	\text{arithmetic theory of } \GL_1.}
\end{eqnarray*}
In Sections \ref{modsym}-\ref{zeta_elts}, we study the map $\varpi$.  In Sections \ref{ordinary} and \ref{mw}, we study the map $\Upsilon$.  In Sections \ref{conj} and \ref{key_ideas}, we state the conjecture and the main result on it.

\subsection{From modular symbols to cup products} \label{modsym}

We construct the map $\varpi_r$ that relates modular symbols in the homology of $X_1(Np^r)$ to cup products in the cohomology of the maximal unramified outside $Np$ extension of $\Q(\zr)$.

\begin{sbpara}

We introduce homology groups $\mcS_r$ and $\tS_r$ of modular curves.

Let $\HH$ denote the complex upper half-plane and $\Gamma_1(Np^r)$ the usual congruence subgroup
of matrices in $\SL_2(\Z)$ that are upper-triangular and unipotent modulo $Np^r$.   We consider the complex
points $Y_1(Np^r) = \Gamma_1(Np^r) \backslash \mb{H}$ of the open modular curve over $\C$.  It is traditional
to use $\cusps$ to denote the cusps $\Gamma_1(Np^r) \backslash \mb{P}^1(\Q)$, but let us instead set
$C_r = \cusps$.  We let 
$$
	X_1(Np^r) = Y_1(Np^r) \cup C_r = \Gamma_1(Np^r) \backslash \HH^*,
$$
 be the closed modular curve, where $\HH^* = \HH \cup \mb{P}^1(\Q)$  is the extended upper half-plane.

The usual modular symbols lie in the first homology group of the space $X_1(Np^r)$ relative to the cusps.  However, $H_1(X_1(Np^r),C_r,\Z)$ is not exactly the natural object for our study.  
Rather, we are interested in its quotient by the action of complex conjugation, the plus quotient.\footnote{This is still not quite the right object unless we invert $2$.  In Section \ref{GL_d}, we take the point of view that the right object is the relative homology of the quotient of the space $X_1(Np^r)$ by the action of complex conjugation.}
We consider the plus quotients of homology and homology relative to the cusps:
\begin{eqnarray*}
	\mcS_r = H_1(X_1(Np^r), \Z)_+ &\mr{and}& \tS_r = H_1(X_1(Np^r), C_r, \Z)_+,
\end{eqnarray*}
where $(\enspace)_+$ denotes the plus quotient.

\end{sbpara}

\begin{sbpara} \label{maninsym}

We introduce Manin symbols $[u:v]_r\in \tS_r$.

Let $u,v\in \Z/Np^r\Z$ be such that $(u,v)=(1)$.
For such $u$ and $v$, we can find $\gamma=\left(\begin{smallmatrix} a & b \\ c & d\end{smallmatrix}\right) \in \SL_2(\Z)$ with $u = c \bmod Np^r$ and $v = d \bmod Np^r$. Define
$$
	[u:v]_r = \left\{ \frac{d}{bNp^r} \to \frac{c}{aNp^r} \right\}_r,
$$
where $\{\alpha \to \beta\}_r$ for $\alpha, \beta \in \mb{P}^1(\Q)$ denotes the class in $\tS_r$ of the hyperbolic geodesic on $\HH^*$ from $\alpha$
to $\beta$.  Then $[u:v]_r$ is independent of the choice of $\gamma$.

By the work of Manin \cite{manin}, we have that the group $\tS_r$ of modular symbols is explicitly presented as an abelian group by generators $[u:v]_r$ and relations
\begin{eqnarray*}
	[u:v]_r= [-u:v]_r = -[v:u]_r  &\mr{and}& [u:v]_r= [u:u+v]_r + [u+v:v]_r.
\end{eqnarray*}

\end{sbpara}

\begin{sbpara} We define an intermediate relative homology group $\tS_r^0$ used in
constructing $\varpi_r$.

We do not use all modular symbols to connect with $\GL_1$.  Rather, we use those modular symbols with boundaries in cusps that do not lie over the cusp at $0$ in $X_0(Np^r) = \Gamma_0(Np^r) \backslash \HH^*$.  Let us denote the set of cusps of $X_1(Np^r)$ that do not lie over the $0$-cusp of $X_0(Np^r)$ by $C_r^0$.  The intermediate space
$$
	\tS_r^0 = H_1(X_1(Np^r), C_r^0, \Z)_+
$$
is the largest space on which we may define $\varpi_r$ and have it factor through the Eisenstein quotient (see \ref{varpi_Eis}).  We have 
$\mcS_r \subset \tS_r^0 \subset \tS_r$.

Our intermediate space also has a simple presentation: it is generated by the $[u:v]$ for nonzero 
$u, v \in \Z/Np^r\Z$ with $(u,v)=(1)$, together with the relations of \ref{maninsym}, again for nonzero $u$ and $v$, and excluding the last relation when $u+v = 0$.

\end{sbpara}

\begin{sbpara}  We define the map $\varpi_r$, which gives our first connection between $\GL_2$ and $\GL_1$.

We start with the primitive $Np^r$th root of unity $\zr = e^{2 \pi i/Np^r}$.
It generates the cyclotomic field $E_r = \Q(\zr)$ and its integer ring $\Z[\zr]$.  Inside 
$E_r$, we have the maximal totally real subfield $F_r = \Q(\zr)^+$ and its integer ring $\mc{O}_r = \Z[\zr]^+$.  

For $a, b \in \Z[\zr,\tfrac{1}{Np}]^{\times}$, we let $\{a,b\}_r$ denote the norm of the Steinberg
symbol of $a$ and $b$ to $K_2(\mc{O}_r[\tfrac{1}{Np}])$.
There is a well-defined homomorphism
$$
	\varpi_r  \colon \tS_r^0 \otimes \Z[\tfrac{1}{2}] \to 
	K_2(\mc{O}_r[ \tfrac{1}{Np}]) \otimes \Z[\tfrac{1}{2}], \qquad
	[u:v]_r \mapsto \{1-\zr^u, 1-\zr^v\}_r
$$
for $u,v \neq 0$.  
Using the Steinberg relation $\{x, 1-x\}_r = 0$ in $K_2$ for $x, 1-x \in \Z[\zr,\tfrac{1}{Np}]^{\times}$,  one
may easily check that the $\{1-\zr^u,1-\zr^v\}_r$ satisfy the same relations as the $[u:v]_r$ 
(see \cite{sh-Lfn, bus} for instance).  It is necessary to invert $2$ for these relations to hold.

\end{sbpara}

\begin{sbpara}\label{kummer}

We interpret $\varpi_r$ on $p$-completions in terms of cup products in Galois cohomology.

For a commutative ring $R$ in which $p$ is invertible, the Kummer exact sequence
$$
	0\to \Z/p^n\Z(1) \to {\mathbb G}_m\overset{p^n}\to {\mathbb G}_m\to 0
$$
on $\Spec(R)_{\et}$ induces the connecting map $R^\times \to H^1_{\et}(R, \Z/p^n\Z(1))$. We have also the Chern class map $K_2(R)\to H^2_{\et}(R, \Z/p^n\Z(2))$.  The value of this map on a product (i.e., Steinberg symbol) in $K_2(R)$ of a pair of elements of $R^{\times}$  is equal to the cup product 
of the images in $H^1_{\et}(R,\Z/p^n\Z(1))$ of the two elements.

We may apply this discussion with $R$ equal to $\Z[\zr,\tfrac{1}{Np}]$ or $\mc{O}_r[\tfrac{1}{Np}]$, 
in which cases the Chern class map 
$K_2(R) \otimes \zp \to H^2_{\et}(R,\zp(2))$
is an isomorphism \cite{tate}.  Moreover, the diagram 
$$
	\SelectTips{eu}{} \xymatrix{ 
	K_2(\Z[\zr,\tfrac{1}{Np}])_+ \otimes \zp \ar[r]^-{\sim} \ar[d]^{\wr}_{N} & 
	H^2_{\et}(\Z[\zr,\tfrac{1}{Np}],\zp(2))_+ \ar[d]^{\wr}_{\cor} \\
	K_2(\mc{O}_r[\tfrac{1}{Np}]) \ar[r]^-{\sim} \otimes \zp & H^2_{\et}(\mc{O}_r[\tfrac{1}{Np}], \zp(2)) }
$$ 
commutes, where $N$ is induced by the norm and $\cor$ is induced by corestriction.  The map $\cor$
is an isomorphism as $\mc{O}_r[\tfrac{1}{Np}]$ has $p$-cohomological dimension $2$.
Let $(1-\zr^u,1-\zr^v)_r$ denote the corestriction of the cup product 
of the elements $1-\zr^u$ and $1-\zr^v$ of
$$
 	\Z[\zr,\tfrac{1}{Np}]^\times \otimes \zp \xsim 
	H^1_{\et}(\Z[\zr,\tfrac{1}{Np}], \Z_p(1)).
$$
By definition of the symbols, the Chern class map in the lower row of the diagram satisfies
$$
	\{1-\zr^u, 1-\zr^v\}_r \mapsto (1-\zr^u,1-\zr^v)_r.
$$ 

We will study the homomorphism to Galois cohomology
$$
	\varpi_r \colon \tS_r^0 \otimes \zp \to H^2_{\et}(\mc{O}_r[\tfrac{1}{Np}],\zp(2)),
	\qquad [u:v]_r \mapsto (1-\zr^u,1-\zr^v)_r,
$$
which is identified with our original $\varpi_r$ on $p$-completions.

\end{sbpara}

\begin{sbpara} \label{Eis}

We define Hecke algebras $\mb{T}_r$ and $\tbT_r$ and their Eisenstein ideals $I_r$ and $\mf{I}_r$.

The Hecke operators $T(n)$ for $n \ge 1$ generate a subalgebra $\tbT_r$ of $\End_{\zp}(\tS_r \otimes \zp)$, the modular Hecke algebra.  We will be interested in this section only in its action on $\tS_r^0$.  We also have a cuspidal Hecke algebra $\mb{T}_r$ of $\End_{\zp}(\mcS_r \otimes \zp)$ and a canonical surjection $\tbT_r  \twoheadrightarrow \mb{T}_r$.  These Hecke algebras contain diamond operators $\langle d \rangle$ for $d \in \Z$, which we take to be $0$ if $(d,Np) \ne 1$.

The Hecke algebra $\tbT_r$ contains the Eisenstein ideal $\mf{I}_r$ generated by the
$T(n) - \sum_{d|n} d\langle d\rangle$ for $n \ge 1$.  It is also generated by $T(\ell) - 1 - \ell \langle \ell \rangle$ for primes $\ell$.  The image $I_r$ of $\mf{I}_r$ in $\mb{T}_r$ is an Eisenstein ideal with the same generators.

\end{sbpara}

\begin{sbpara} \label{varpi_Eis}

We connect our study of $\varpi_r$ with the Eisenstein ideal.

The third author conjectured \cite{sh-Lfn} (on $\mcS_r$, see also \cite{bus} for $N = 1$), and the first two authors proved \cite[Theorem 5.2.3]{fk-proof} that
$\varpi_r$ is ``Eisenstein.''  By this, we mean
that $\varpi_r$ factors through the quotient of $\tS_r^0 \otimes \zp$ by the Eisenstein ideal, that is, through a map
$$
	(\tS_r^0/\mf{I}_r\tS_r^0) \otimes \zp \to H^2_{\et}(\mc{O}_r[\tfrac{1}{Np}],\zp(2)).
$$
We can show that  this follows from the fact that $\varpi_r$ is the specialization of a map in the 
$\GL_2$-setting: see Section \ref{zeta_elts}.

\end{sbpara}

\begin{sbpara}

Let $G_r= (\Z/Np^r\Z)^\times/\{\pm 1\}$, and set $\Lambda_r = \zp[G_r]$.  The algebra $\Lambda_r$ appears in 
two different contexts in our story:

\begin{enumerate}
	\item[(1)] On the $\GL_2$-side, $\Lambda_r$ is a $\zp$-algebra of diamond operators in 
	$\mb{T}_r$ (or $\tbT_r$):
	we define a $\zp$-linear injection 
	$\iota_r \colon \Lambda_r \hookrightarrow \mb{T}_r$
	that sends the group element in $\Lambda_r$ corresponding to $a \in (\Z/Np^r\Z)^{\times}/\langle -1 \rangle$ 
	to the \underline{inverse} 
	$\langle a \rangle^{-1}$ of the diamond operator for $a$ (i.e., for any lift of $a$ to an integer).
	
	\item[(2)] On the $\GL_1$-side, $\Lambda_r$ is the $\zp$-group ring of $\Gal(F_r/\Q)$: 
	we have an isomorphism
	$$
		(\Z/Np^r\Z)^\times \xsim \Gal(E_r/\Q), \qquad a \mapsto \sigma_a, 
	$$ 
	where $\sigma_a(\zr)=\zr^a$.  This
	gives rise to an isomorphism $G_r  \xsim \Gal(F_r/\Q)$ that is the map on group elements
	defining $\Lambda_r \xsim \zp[\Gal(F_r/\Q)]$.
\end{enumerate}

These actions are compatible with $\varpi_r$ in the sense that for any $x \in \tS_r^0 \otimes \zp$ and
$a \in (\Z/Np^r\Z)^{\times}$, we have 
$$
	\varpi_r(\langle a \rangle^{-1} x) = \sigma_a \varpi_r(x).
$$
This is easily seen: taking $x = [u:v]_r$ for some nonzero $u$ and $v$, we have
\begin{eqnarray*}
	\langle a \rangle^{-1}[u:v]_r = [au:av]_r  &\mr{and}& \sigma_a(1-\zr^u,1-\zr^v)_r 
	= (1-\zr^{au},1-\zr^{av})_r.
\end{eqnarray*}
So, to say that $\varpi_r$ is Eisenstein is to say that $\varpi_r(T(\ell)x) = (1+\ell\sigma_{\ell}^{-1})\varpi_r(x)$
for primes $\ell \nmid Np$ and $\varpi_r(T(\ell)x) = \varpi_r(x)$ for $\ell \mid Np$.

\end{sbpara}

\subsection{Passing up the modular and cyclotomic towers: the map $\varpi$} \label{tower}

We pass up the modular 
tower on the $\GL_2$-side and the cyclotomic tower on the $\GL_1$-side to define the map $\varpi = \varprojlim_r\varpi_r$.  

\begin{sbpara} \label{Lambda}

Let $G = \varprojlim_r G_r$.  Then the completed group ring
$$
	\La = \zp\ps{G} = \vpr \Lambda_r
$$
is the Iwasawa algebra for $G$.  As with $\Lambda_r$, let us emphasize its dual nature:
\begin{enumerate}
	\item[(1)] Set $\mb{T} = \varprojlim_r \mb{T}_r$ and $\tbT = \varprojlim_r \tbT_r$.  
	The projective limit of the injections $\iota_r$
	defines a map $\iota \colon \Lambda \hookrightarrow \mb{T}$ of profinite $\zp$-modules that takes
	$a \in G$ to the projective system of inverses $\langle a \rangle^{-1}$ of diamond operators 
	corresponding to $a$.
	\item[(2)] Set $K = \cup_{r \ge 1} F_r$, the maximal totally real subfield of $L = \cup_{r \ge 1} E_r$.  Then our
	identifications $\Lambda_r \xsim \zp[\Gal(F_r/\Q)]$ for $r \ge 1$ induce an isomorphism
	$\Lambda \xsim \zp\ps{\Gal(K/\Q)}$ of completed group rings in the projective limit.
\end{enumerate}

\end{sbpara}

\begin{sbpara}\label{limitpi}

We have the following projective limits of spaces of modular symbols:
\begin{eqnarray*}
	\mcS = \vpr\,(\mcS_r \otimes \zp) &\mr{and}& \tS^0 = \vpr\,(\tS_r^0 \otimes \zp).
\end{eqnarray*}
Let $\mf{I} \subset \tbT$ and $I \subset \mb{T}$ be the Eisenstein ideals, defined by the same set of generators as $I_r$, but now viewed as compatible sequences of operators in the Hecke algebras.

Our maps $\varpi_r$ are compatible with change of $r$ and induce in the projective limit a map
$$
	\varpi \colon \tS^0 \to \vpr H^2_{\et}(\mc{O}_r[\tfrac{1}{Np}], \zp(2))
$$
that factors through $\tS^0/\mf{I}\tS^0$ by the result of \cite{fk-proof}.  This map $\varpi$ is a homomorphism of $\Lambda$-modules, the actions arising from part (1) of \ref{Lambda} on the left and part (2) of \ref{Lambda} on the right.

\end{sbpara}

\begin{sbpara} \label{XY}

We recall the unramified Iwasawa module $X$, study its difference from Galois cohomology, and consider a related $\Lambda$-module $Y$.

Let $X$ be the projective limit of the $p$-parts $A_r$ of the ideal class groups of the fields $E_r$.  Class field theory allows us to identify $X$ with the Galois group of the maximal unramified abelian pro-$p$ extension of $\tK$.

For $R$ as in \ref{kummer}, the Kummer exact sequence induces 
$$
	\Pic(R)=H^1_{\et}(R, {\mathbb G}_m)\to H^2_{\et}(R, \Z/p^n\Z(1)).
$$ 
Taking a projective limit of such maps for the rings $R=\mc{O}_r[\tfrac{1}{Np}]$, we obtain
$$
	X = \vpr A_r \to \vpr H^2_{\et}(\mc{O}_r[ \tfrac{1}{Np}], \Z_p(1)).
$$ 
In general, this map is neither injective nor surjective.  Its kernel and cokernel can be explicitly described as contributions of classes of primes and Brauer groups at places dividing $Np$, respectively.  
We will deal with a part of cohomology on which this subtle difference disappears.

The Iwasawa algebra $\zp\ps{\Gal(\tK/\Q)}$ acts on $X$, but this action does not in general factor through $\Lambda$.  We want to consider the $(-1)$-eigenspace $X^-$ of $X$ under complex conjugation.  To do so, we take the Tate twist $Y = X^-(1)$, or equivalently, the fixed part $X(1)^+$.  Then $\sigma_{-1}$ acts trivially on $Y$, so $Y$ is a $\La$-module.  The map from $X$ to cohomology induces a $\Lambda$-module homomorphism 
$$
	 Y \to \vpr H^2_{\et}(\mc{O}_r[\tfrac{1}{Np}], \zp(2)).
$$
Together with $\varpi$, this will allow us to relate $\mcS$ with $Y$.

\end{sbpara}

\begin{sbpara} \label{varpi}

We have two objects of study:
\begin{enumerate}
	\item[$\bullet$] the geometric object $P = \mcS/I\mcS$ for $\GL_2$,
	\item[$\bullet$] the Iwasawa-theoretic object $Y = X^-(1)$ for $\GL_1$.
\end{enumerate}
We can relate these on $\theta$-parts for suitable even characters $\theta$ of $(\Z/Np\Z)^{\times}$.

For a primitive, even character $\theta \colon (\Z/Np\Z)^{\times} \to \qpbt$, we may consider the quotient $\Lambda_{\theta} = \Lambda \otimes_{\zp[(\Z/Np\Z)^{\times}]} \zp[\theta]$ of $\La$, where $\zp[\theta]$ is the $\zp$-algebra generated by the values of $\theta$.  For a $\Lambda$-module $M$, we then
let $M_{\theta} = M \otimes_{\Lambda} \Lambda_{\theta}$ denote its $\theta$-part.

We need a technical assumption to insure that the maps
\begin{eqnarray*}
	P_{\theta} \to \tS^0_{\theta}/\mf{I}_{\theta}\tS^0_{\theta} &\mr{and}& 
	Y_{\theta} \to \vpr H^2_{\et}(\mc{O}_r[\tfrac{1}{Np}], \zp(2))_{\theta}
\end{eqnarray*}
are isomorphisms.  Together with primitivity, the assumption is as follows:
\begin{enumerate}
	\item[$\bullet$] $\theta\omega^{-1}|_{(\Z/p\Z)^{\times}} \neq 1$ or 
	$\theta\omega^{-1}|_{(\Z/N\Z)^{\times}}(p) \neq 1$,
\end{enumerate}  
where $\omega \colon 
(\Z/Np\Z)^{\times} \to \zp^{\times}$ denotes the Teichm\"uller character (i.e., projection to 
$(\Z/p\Z)^{\times}\subset \zp^{\times}$).
For such a $\theta$, our $\varpi$ induces a map on $\theta$-parts 
$\varpi \colon \mcS_{\theta} \to Y_{\theta}$ that will factor through $P_{\theta}$.

\end{sbpara}

\subsection{Zeta elements: $\varpi$ is ``Eisenstein''} \label{zeta_elts}

We sketch the proof that $\varpi$ factors through the quotient of $\tS^0$ by the 
Eisenstein ideal $\mf{I}$.

\begin{sbpara}
Let $Y_1(Np^r)$ be the moduli space of pairs $(E, e)$ where $E$ is an elliptic curve and $e$ is a point of order $Np^r$ on $E$, and let $Y(Np^r)$ be the moduli space of elliptic curves endowed with a full $Np^r$-level structure.  We view these moduli spaces as schemes over $\Z[\tfrac{1}{Np}]$.
For any nonzero $(\alpha, \beta) \in \tfrac{1}{Np^r}\Z^2/\Z^2$, there is a Siegel unit
$g_{\alpha,\beta} \in \mc{O}(Y(Np^r))^{\times}$.\footnote{Actually, $g_{\alpha,\beta}$ is a root of 
a unit, but the difficulties this causes are resolvable by passing to the projective limit and descending, so 
we ignore this for simplicity of presentation.  We will be very careless about
denominators in several places, omitting them where they occur for simplicity of the discussion that follows.}  
It has the $q$-expansion
$$
	g_{\alpha,\beta} = q^{\frac{1}{12}-\frac{\alpha}{2}+\frac{\alpha^2}{2}}
	\prod_{n=0}^{\infty}(1-q^{n+\alpha}e^{2\pi i \beta})
	\prod_{n=1}^{\infty}(1-q^{n-\alpha}e^{-2\pi i \beta}) \in 
	\mc{O}_r[\tfrac{1}{Np}]\ps{q^{1/12Np^r}}[q^{-1}]^{\times}.
$$
If $\alpha = 0$, then we may view $g_{0,\beta}$ as an element of $\mc{O}(Y_1(Np^r))^{\times}$.
The crucial point is that the specialization of a Siegel unit
of the form $g_{0,\frac{u}{Np^r}}$ at the $\infty$-cusp is the cyclotomic $Np$-unit $1-\zr^u$.
Specifically, this specialization is given by projecting its $q$-expansion to 
$\mc{O}_r[\frac{1}{Np}]\ps{q^{1/Np^r}}^{\times}$ and then evaluating at $q = 0$ \cite[Section 5.1]{fk-proof}.

\end{sbpara}

\begin{sbpara} \label{z_r}
We have a homomorphism
$$
	z_r \colon \tS_r^0 
	\to H^2_{\et}(Y_1(Np^r), \zp(2)), \qquad
	z_r([u:v]_r) = g_{0,\frac{u}{Np^r}} \cup g_{0,\frac{v}{Np^r}}
$$
of $\tbT_r$-modules that takes a Manin symbol to a Beilinson element given by a cup product of two Siegel units 
\cite[Proposition 3.3.15]{fk-proof}.  Related elements were studied in \cite{kato}.
\end{sbpara}

\begin{sbpara} \label{infty_r}

There is again a specialization-at-$\infty$ map 
$$
	\infty_r \colon H^2_{\et}(Y_1(Np^r), \zp(2)) \to H^2_{\et}(\mc{O}_r[ \tfrac{1}{Np}],\zp(2))
$$
that takes $g_{0,\frac{u}{Np^r}} \cup g_{0,\frac{v}{Np^r}}$ to 
$(1-\zr^u,1-\zr^v)_r$.  
So, cup products of cyclotomic units are specializations at cusps of Beilinson elements. We have
$$\varpi_r = \infty_r \circ z_r \colon \tS_r^0 \to H^2_{\et}(\mc{O}_r[ \tfrac{1}{Np}], \Z_p(2)).$$
It can be shown that specialization at $\infty$ is Eisenstein. Hence so is 
$\varpi_r$ \cite[Sections 5.1-5.2]{fk-proof}.  

\end{sbpara}

\begin{sbpara}

By passing the projective limit over $r$, we see that $\varpi$ is Eisenstein.  The identity $\varpi = \infty \circ z$ is the commutativity of the left-hand square in the diagram of \ref{key_diag}.

\end{sbpara}

\subsection{Ordinary homology groups of modular curves} \label{ordinary}

Homology groups of the modular curves are useful for us in two different ways.  They contain modular symbols, allowing us to define $\varpi$.  They also have Galois actions, allowing us to define $\Upsilon$, which is our next goal. We use two 
different groups derived from homology, $\mcS$ as above and $\mcT$ defined below, to construct the two maps.
For the modular symbols, we require only the plus part of homology.
On the other hand, to have Galois actions, we cannot restrict to plus parts.  Instead,
we take ordinary parts to control the growth of homology groups in the modular tower 
and to specify the form of the local Galois action at $p$.  
The fact that we use different groups should be kept in mind in the $\GL_d$-setting, in which we will not consider $\Upsilon$.

\begin{sbpara}

We introduce Hida's ordinary $p$-adic cuspidal and modular Hecke algebras $\mf{h}$ and $\mf{H}$. 

Recall our cuspidal Hecke algebra $\mb{T}$ from \ref{Lambda}, which acts $\Lambda$-linearly on $\mcS$.  The action of $T(p)$ breaks it into a direct product of two rings: an ordinary part in which the image of $T(p)$ is invertible and another part in which $T(p)$ is topologically nilpotent.  The ordinary cuspidal $p$-adic Hecke algebra $\mf{h} = \mb{T}^{\ord}$ of Hida \cite{hida-iwa} is this ordinary part.  This is a $\Lambda$-subalgebra that is projective of finite $\Lambda$-rank.  We may speak of Hecke operators $T(n) \in \mf{h}$ by taking the images
of the $T(n) \in \mb{T}$.

The Hecke algebra $\mf{h}$ is remarkable in that it simply encapsulates information about the ordinary Hecke algebras of all weights $\ge 2$ and all levels dividing some $Np^r$.  For instance, its quotient for the action of the kernel of $G \to G_r$ is the Hecke algebra $\mf{h}_r = \mb{T}_r^{\ord}$.  This highly regular behavior is the subject of Hida theory.

We also have the ordinary modular Hecke algebra $\mf{H} = \tbT^{\ord}$, of which $\mf{h}$ is a quotient.  In general, if $M$ is a $\tbT$-module (resp., $\mb{T}$-module), then we use $M^{\ord}$ to denote its ordinary part, the maximal summand on which $T(p)$ acts invertibly, which is an
$\mf{H}$-module (resp., $\mf{h}$-module).

\end{sbpara}

\begin{sbpara} \label{homology-hecke} 

We introduce the ordinary homology groups $\mcT$ and $\tT$.  These have commuting actions of Hecke algebras and the absolute Galois group $G_{\Q} = \Gal(\qbar/\Q)$.  The study of these actions on $\mcT$ will allow us to define the map $\Upsilon$ in Section \ref{mw}.

The Hecke operators $T(n)$ with $n \ge 1$ act on the homology of $X_1(Np^r)(\C)$ and the homology relative to the cusps and are compatible with projective limits.
We consider the ordinary parts $\mcT$ and $\tT$ of the projective limits
$$
	\mcT = \vpr H_1(X_1(Np^r), \Z_p)^{\ord} \subset
	\tT = \vpr H_1(X_1(Np^r), C_r, \Z_p)^{\ord}.
$$
We are primarily interested in $\mcT$.  The $\mbT$-action on $\mcT$ factors through $\mf{h}$.  As an $\mf{h}$-module, $\mcT$ is finitely generated and torsion-free, and $\mcT$ is projective of finite rank over $\La$.  If we denote by $Q(\La)$ the total quotient ring of $\La$, then $Q(\La)\otimes_{\La}\mcT $ is a free $Q(\La) \otimes_{\La} \mf{h}$-module of rank $2$. 

The absolute Galois group $G_\Q= \Gal(\qbar/\Q)$ acts on the homology of $X_1(Np^r)$ by its duality with cohomology and the identification of Betti cohomology with \'etale cohomology of the scheme $X_1(Np^r)_{/\qbar}$ over $\qbar$.  This describes the first of the two isomorphisms
\begin{eqnarray*} 
	&H_1(X_1(Np^r), \Z_p)\cong \text{Hom}(H^1_{\et}(X_1(Np^r)_{/\qbar}, \Z_p), \Z_p),&\\
	&H_1(X_1(Np^r), C_r, \Z_p)\cong \text{Hom}(H^1_{\et,c}(Y_1(Np^r)_{/\qbar}, \Z_p), \Z_p),&
\end{eqnarray*}
where in the second, the duality of the relative cohomology group is with the compactly supported cohomology of the open modular curve.  This Galois action commutes with the action of the Hecke operators, so passes to ordinary parts, and it is compatible in the towers.  Therefore, $\mf{H}[G_{\Q}]$ acts compatibly on $\mcT$ and $\tT$.

\end{sbpara}

\begin{sbpara}

We introduce the Eisenstein ideals $I$ and $\mf{I}$ of the ordinary Hecke algebras.  
 
Let us reuse the notation $I$, allowing it to denote the Eisenstein ideal of $\mf{h}$, which is the image of the Eisenstein ideal $I$ of $\mb{T}$ in $\mf{h}$.  
We remark that, since $T(p)-1 \in I$ and $1$ is a unit, the quotient map $\mb{T}/I \to \mf{h}/I$ is an isomorphism. 
We will also reuse the notation $\mf{I}$ for the Eisenstein ideal of $\mf{H}$, the image of $\mf{I} \subset \tbT$. 
 
\end{sbpara}

\begin{sbpara} \label{two_primes}

In the $\GL_2$-setting over $\Q$, there are two places which play important roles:  the place 
at $p$ and the real place.  We study the actions of the corresponding local Galois groups.

We first study the local action at $p$: here we have an interesting quotient $\mcT_{\quo}$.  The fact that $\mcT$ is ordinary for $T(p)$ tells us about the action of $G_{\qp}$, which is to say that it is ordinary in the sense of $p$-adic Hodge theory.  More specifically to our case, we have an exact sequence
$$
	0 \to \mcT_{\sub} \to \mcT \to \mcT_{\quo} \to 0
$$
of $\mf{h}[G_{\qp}]$-modules, with $\mcT_{\sub}$ and $\mcT_{\quo}$ defined as follows.
First, $\mcT_{\sub}$ is the largest submodule of $\mcT$ such that $G_{\qp}$ acts on $\mcT_{\sub}(-1)$ by inverse diamond operators, and $\mcT_{\quo}$ is the quotient.   Put more simply,
$\mcT_{\quo}$ is the maximal unramified, $\mf{h}$-torsion-free quotient of $\mcT$.  

At the real place, we have $\mcT^+$, which is isomorphic to $\mcS^{\ord}$.  It fits in an exact sequence
$$
	0 \to \mcT^+ \to \mcT \to \mcT/\mcT^+ \to 0
$$
of $\mf{h}[G_{\R}]$-modules, and both $Q(\La) \otimes_{\La} \mcT^+$ and $Q(\La) \otimes_{\La} \mcT/\mcT^+$ 
are free of rank $1$ over $Q(\La) \otimes_{\La} \mf{h}$. 

The  compositions $\mcT^+ \to \mcT \to \mcT_{\quo}$ and $\mcT_{\sub} \to \mcT \to \mcT/\mcT^+$ relate 
the two exact sequences.  We study these maps on Eisenstein components in \ref{unramified}.    The interplay between the reductions modulo $I$ of the two exact sequences allows us to construct the map $\Upsilon$.

\end{sbpara}

\begin{sbpara} \label{Lambda-adic}

We discuss $\La$-adic cusp forms and modular forms and their ordinary parts $\mf{S}$ and 
$\mf{M}$.  

Let $S_2(Np^r)_{\Z}$ denote the space of cusp forms of weight $2$ and level $Np^r$
with integer coefficients.  For a ring $R$, we then set $S_2(Np^r)_R = S_2(Np^r)_{\Z} \otimes R$.  If $\epsilon \colon G_r \to R^{\times}$ is a homomorphism, then we may speak of $S_2(Np^r,\epsilon)_R$, those cusp forms
in $S_2(Np^r)_R$ with nebentypus $\epsilon$.

Any finite order character $\epsilon \colon G \to \qpbt$ induces a ring homomorphism $\La \to \qpbar$.  We let $\tilde{\epsilon} \colon \La\ps{q} \to \qpbar\ps{q}$ be the induced map on coefficients.  
An element $f \in \La\ps{q}$ is said to be a $\La$-adic cusp form of weight $2$ and level $Np^{\infty}$ if for every $\epsilon$, one has
$\tilde{\epsilon}(f) \in S_2(Np^r,\epsilon)_{\qpbar}$ 
with $r \ge 0$ such that $\epsilon$ factors through $G_r$ \cite{wiles, ohta-eich}.
We denote the set of such $\La$-adic cusp forms by $S_{\La}$.

The Hecke operators $T(n)$ for $n \ge 1$ act on $S_{\La}$ via the usual formal action of Hecke operators on $q$-expansions.  We define $\mf{S}$ to be the ordinary part $S_{\La}^{\mr{ord}}$ of $S_{\La}$.\footnote{There is one potentially confusing aspect: the action of $\Lambda \hookrightarrow \mf{h}$ on 
$\mf{S} \subset \La\ps{q}$ is not given by multiplication of the coefficients of $q$-expansions by the element 
of $\Lambda$.  It is instead this multiplication after first applying the inversion map $\lambda \mapsto \lambda^*$ on $\La$ that takes group elements to their inverses.}
The ordinary $\La$-adic cusp forms and the ordinary Hecke algebra are dual in the usual sense.  That is, we have a perfect pairing of $\La$-modules,
$$
	\mf{h} \times \mf{S} \to \La, \qquad (T,f) \mapsto a_1(Tf),
$$
where $a_1(g)$ denotes the $q$-coefficient in the $q$-expansion of $g \in S_{\La}$.
As a consequence, $Q(\Lambda) \otimes_{\La} \mf{S}$ is free of rank $1$ over $Q(\Lambda) \otimes_{\La} \mf{h}$.

Similarly, we have a space $\mf{M}$ of ordinary $\La$-adic modular forms with $q$-expansions that are integral outside of the constant term, which sits inside $Q(\Lambda) + \La\ps{q}$.  There is a perfect pairing
$\mf{H} \times \mf{M} \to \La$ that restricts to the pairing for cusp forms.

\end{sbpara}

\begin{sbpara} \label{ohta}

As we shall explain in a more canonical fashion in \ref{ES}, there is an isomorphism
$\mcT_{\quo} \cong \mf{S}$ of $\mf{h}$-modules given by Ohta's $\La$-adic Eichler-Shimura isomorphism \cite{ohta-eich, ohta-ord2}.   Moreover, Ohta showed that $\mcT_{\sub} \cong \mf{h}$ via a $\La$-duality with $\mcT_{\quo}$.

\end{sbpara}

\subsection{Refining the method of Ribet and Mazur-Wiles: the map $\Upsilon$} \label{mw}

We define the map $\Upsilon$ of \cite{sh-Lfn} 
and consider the relationship with the work of Mazur-Wiles \cite{mw}.  
Our description is heavily influenced by the approaches of Wiles \cite{wiles} and Ohta \cite{ohta-ord2}.

We suppose that $p \ge 5$ and $p \nmid \varphi(N)$.\footnote{It should actually be possible to allow either or both of $p = 3$ and $p \mid \varphi(N)$ in what follows.}
We will work mostly in the $\theta$-part (as in \ref{varpi}) for a fixed primitive, even character $\theta \colon (\Z/Np\Z)^{\times} \to \qpbt$ such that the condition $\theta\omega^{-1}|_{(\Z/p\Z)^{\times}} \neq 1$ or $\theta\omega^{-1}|_{(\Z/N\Z)^{\times}}(p) \neq 1$ of \ref{varpi} holds.  We also suppose that $\theta \neq \omega^2$ in the case that $N = 1$.

\begin{sbpara} \label{Ups_idea}

We briefly outline the construction of $\Upsilon \colon Y_{\theta} \to P_{\theta}$ that will appear in this section.

We analyze the $\mf{h}[G_{\Q}]$-action on $\mcT_{\theta}/I_{\theta}\mcT_{\theta}$, showing that it fits in an exact sequence
$$
	0 \to P_{\theta} \to \mcT_{\theta}/I_{\theta}\mcT_{\theta} \to Q_{\theta} \to 0
$$
of $\mf{h}[G_{\Q}]$-modules.  Any such exact sequence provides a cocycle $G_{\Q} \to \Hom_{\mf{h}}(Q_{\theta},P_{\theta})$ that defines its extension class in Galois cohomology.  Our exact sequence has three key properties: the $G_{\tK}$-action on $\mcT_{\theta}/I_{\theta}\mcT_{\theta}$ is unramified, the $G_{\tK}$-actions on $P_{\theta}$ and $Q_{\theta}$ are trivial, and the $\mf{h}_{\theta}/I_{\theta}$-module $Q_{\theta}$ is free of rank $1$ with a canonical generator.  We may therefore modify our cocycle as follows.  First, we compose it with evaluation at the generator of $Q_{\theta}$ to obtain a map $G_{\Q} \to P_{\theta}$.  Since $G_{\tK}$ acts trivially on $P_{\theta}$ and $Q_{\theta}$, this map restricts to a homomorphism $G_{\tK} \to P_{\theta}$.  Since the $G_{\tK}$-action on $\mcT_{\theta}/I_{\theta}\mcT_{\theta}$ is unramified, this homomorphism in turn factors through a homomorphism $X \to P_{\theta}$.  After a twist, it further factors through $Y_{\theta}$ and provides the desired map $\Upsilon \colon Y_{\theta} \to P_{\theta}$, which we can show to be of $\La$-modules.  

We first explain that $\mf{h}_{\theta}/I_{\theta}$ is the quotient of $\Lambda_{\theta}$ by a $p$-adic $L$-function $\xi_{\theta}$.  This will provide the connection between the map $\Upsilon$ and the Iwasawa main conjecture.

\end{sbpara}

\begin{sbpara} \label{xi}

We define the $p$-adic $L$-function $\xi_{\theta}$.

Note that any homomorphism $G \to \qbt$ factors through some $G_r$ and so induces an even Dirichlet character.  Note also that $G_1 = (\Z/Np\Z)^{\times}/\langle -1 \rangle$ and $G \cong G_1 \times (1+p\zp)$.

The $p$-adic $L$-function $\xi_{\theta}$ is the unique element of $\La_{\theta}$ that interpolates Dirichlet $L$-values at $-1$ in the sense that for each character $\epsilon \colon G \to \qbt$ such that $\epsilon|_{G_1} = \theta$, the ring homomorphism $\Lambda_{\theta} \to \qpbar$ induced by $\epsilon$ sends $\xi_{\theta}$ to the value $L(\epsilon^{-1},-1) \in \qbar$ of the Dirichlet $L$-function.

We can also describe $\xi_{\theta}$ in terms of Kubota-Leopoldt $p$-adic $L$-functions.  
We make the identification $\Lambda_{\theta} = \zp[\theta]\ps{T}$ with $T = [u]-1$,
where $[u]$ is the group element of $u \in 1+p\zp$ with $p$-adic logarithm $(1-p^{-1})^{-1}$.  We then have the following equality of functions of $s \in \zp$:
$$
	\xi_{\theta}(u^s-1) = L_p(\omega^2\theta^{-1},s-1).
$$

\end{sbpara}

\begin{sbpara} \label{hmodI}

We construct a canonical isomorphism 
$\mf{h}_{\theta}/I_{\theta} \xsim \Lambda_{\theta}/(\xi_{\theta})$.

Consider the ordinary $\La$-adic Eisenstein series
$$
	\mc{E}_{\theta} = \tfrac{1}{2}(\xi_{\theta})^* + \sum_{n=1}^{\infty} \Bigg( 
	\sum_{\substack{d \mid n \\ (d,Np) = 1}}
	d [d] \Bigg) q^n \in \mf{M}_{\theta},
$$
where $[d]$ is the image in $\Lambda_{\theta^{-1}}$ of the group element in $G$ for $d$, and $\lambda \mapsto \lambda^*$ is the involution defined in the footnote of \ref{Lambda-adic}.\footnote{The reader may wish to ignore
the involutions in order to focus on the idea of the argument.}
By duality with the Hecke algebra, it provides a surjective homomorphism $\mf{H}_{\theta} \to \Lambda_{\theta}$, the kernel of which is $\mf{I}_{\theta}$ by definition.

Let $M_{\Eis}$ denote the component of an $\mf{H}$-module 
$M$ for the unique maximal ideal $\mf{m}$ containing the Eisenstein ideal $\mf{I}_{\theta}$.
By our choice of $\theta$, the Eisenstein series $\mc{E}_{\theta}$ is not congruent modulo $\mf{m}$ to
any other Eisenstein series \cite[Lemma 1.4.9]{ohta-cong}.  
It follows from this that the injection of $\mf{S}$ in $\mf{M}$ induces an exact sequence
$$
	0 \to \mf{S}_{\Eis} \to \mf{M}_{\Eis} \to \Lambda_{\theta} \to 0,
$$
where the latter map takes a modular form to the (involution of the) constant term in its $q$-expansion.  Our map
$\pi_{\theta} \colon \mf{h}_{\theta}/I_{\theta} \to \Lambda_{\theta}/(\xi_{\theta})$ may then be constructed from the reduction of $\mc{E}_{\theta}$ modulo $\xi_{\theta}$.  That is, $\mc{E}_{\theta}$ is a cusp form modulo $(\xi_{\theta}) \subseteq \Lambda_{\theta}$ by the exact sequence, and this cusp form provides the surjective map $\pi_{\theta}$ by duality with the Hecke algebra $\mf{h}$.  Once we know that $\pi_{\theta}$ is an isomorphism, it is inverse to the map induced by $\iota_{\theta}$, where $\iota_{\theta}$ is the $\theta$-part of the map $\iota$ defined in \ref{Lambda}.

We explain the idea behind the injectivity of $\pi_{\theta}$.  We have an evident surjection $\Lambda_{\theta} \to \mf{h}_{\theta}/I_{\theta}$ given by the fact that every Hecke operator $T(n)$ is identified modulo $I_{\theta}$ with an element of $\Lambda_{\theta}$.  So, $\mf{h}_{\theta}/I_{\theta}$ is some quotient of $\Lambda_{\theta}$.  The $\La$-adic forms in $\mf{M}_{\theta}$ have integral constant coefficients, which can be seen by the method of \cite[Proposition 1]{emerton}.  Given this, the existence of $\pi_{\theta}$ is equivalent to the fact $\mc{E}_{\theta}$ modulo $(\xi_{\theta})$ is a $\Lambda$-adic cusp form.  
As the constant coefficient of $\mc{E}_{\theta}$ equals $\xi_{\theta}^*$ times a unit, no surjective homomorphism to a larger quotient of $\Lambda_{\theta}$ can exist.\footnote{Another, more usual, way to approach injectivity is to use $I_{\theta} + \xi_{\theta}\mf{h}_{\theta}$ in place of $I_{\theta}$ until one recovers the equality of these ideals through a proof of the main conjecture, as in \ref{MW1} below.}    Thus, $\pi_{\theta}$ is an isomorphism.

\end{sbpara}

\begin{sbpara} \label{Q}

We define $Q_{\theta}$ and construct a canonical surjection $\mcT_{\theta}/I_{\theta}\mcT_{\theta} \to Q_{\theta}$ of $\mf{h}[G_{\Q}]$-modules.

For a module $M$ over a $\La$-algebra $h$, let $M^{\,\sharp}$ denote the $h[G_{\Q}]$-module that is $M$
as an $h$-module and on which $\sigma \in G_{\Q}$ acts through multiplication by the inverse of the image of $\sigma$ in $G$.
We then define $Q_{\theta} = (\mf{h}_{\theta}/I_{\theta})^{\,\sharp}(1)$.  
Consider the Jacobian variety $J_r$ of the curve $X_1(Np^r)$.  Let $J_{r,\tor} \subset
J_r(\qbar)$ be its torsion subgroup, and take the contravariant 
(i.e., dual)
action of $\mb{T}_r$ on $J_{r,\tor}$.
Consider the class $\alpha_r \in J_r(\C)$  of the divisor $(0)-(\infty)$, where $0$ and $\infty$ are viewed as cusps on $X_1(Np^r)(\C)$.  It is torsion by the theorem of Drinfeld and Manin \cite{drinfeld, manin}.
Moreover, $\alpha_r$ is easily seen to be annihilated by $I_r$.  

Let $\beta_{r,\theta}$ be the image of $\alpha_r$ in the $\theta$-part of $J_r[p^{\infty}] = J_{r,\tor} \otimes
\zp$.  The $\mb{T}_{r,\theta}$-span $B_{r,\theta}$ of $\beta_{r,\theta}$ is a quotient of $\mf{h}_{r,\theta}/I_{r,\theta}$ by definition.  Moreover, $B_{r,\theta}$ is isomorphic to $\La_{r,\theta}/(\xi_{r,\theta})$ by a computation of divisors of Siegel units that says in particular that the $\theta$-part of the divisor of $g_{0,\frac{1}{Np^r}}$ is $\xi_{r,\theta}$ times $(0)-(\infty)$, up to a unit (see \cite[Section 4.2]{mw}).\footnote{In the projective limit, this gives another way of defining the isomorphism $\mf{h}_{\theta}/I_{\theta} \xsim \Lambda_{\theta}/(\xi_{\theta})$.}  
Here, $\xi_{r,\theta}$ denotes the image of $\xi_{\theta}$ in $\Lambda_{r,\theta}$.
The $G_{\Q}$-action on $B_{r,\theta}$ factors through $\Gal(F_r/\Q)$, and we have 
$\sigma_a\beta_{r,\theta} = \langle a \rangle^{-1} \beta_{r,\theta}$ 
for any $a \in (\Z/Np^r\Z)^{\times}$.

Poincar\'e duality allows us to identify the first \'etale homology group of $X_1(Np^r)_{/\qbar}$ with
the Tate twist of the first \'etale cohomology group.  Taking this together with the canonical pairing of
cohomology and the torsion in $J_r$, we obtain a Galois-equivariant, perfect pairing
$$
	(\enspace,\;\;) \colon H_1^{\et}(X_1(Np^r)_{/\qbar}, \Z_p) \times J_r[p^{\infty}] \to
	\qp/\zp(1)
$$ 
with respect to which the Hecke operators are self-adjoint.  Let $(\enspace,\;\,)_{\theta}$ denote the induced pairing on $\theta$-parts.
Define a map $\phi$ by
$$
	\phi \colon H_1^{\et}(X_1(Np^r)_{/\qbar}, \Z_p)_{\theta} \to \Lambda_{r,\theta} \otimes \qp/\zp(1), \qquad 
	x\mapsto \sum_{a\in G_r} [a]_r \otimes (x, \langle a \rangle \beta_{r,\theta})_{\theta},
$$ 
where $[a]_r \in \Lambda_{r,\theta}$ denotes the group element for $a$.\footnote{To make sense of this, note that the tensor product in the sum is taken over $\zp[\theta]$.}  Let $\xi_{r,\theta}$ be the image of $\xi_{\theta}$ in 
$\La_{r,\theta}$.  As $\xi_{r,\theta}\beta_{r,\theta} = 0$, the image of the map $\phi$ is contained in the group $(\Lambda_{r,\theta} \otimes \qp/\zp(1))[\xi_{r,\theta}]$ of $\xi_{r,\theta}$-torsion, and $\phi$
factors through the quotient $\mcT_{r,\theta}/I_{r,\theta}\mcT_{r,\theta}$.

Consider the composition
$$
	\mcT_{r,\theta}/I_{r,\theta}\mcT_{r,\theta} \xrightarrow{\phi} 
	(\Lambda_{r,\theta} \otimes \qp/\zp(1))[\xi_{r,\theta}] \to 
	(\La_{r,\theta}/\xi_{r,\theta})(1) \xrightarrow{\iota_{r,\theta}} (\mf{h}_{r,\theta}/I_{r,\theta})(1),
$$ 
where the second map is given by
$x \mapsto \xi_{r,\theta} \tilde{x}$ for any lifting $\tilde{x}$ of $x$ to $\Q_p[G_r]_{\theta}(1)$.
It is surjective by our description of $B_{r,\theta}$ and the perfectness of $(\enspace,\;\,)_{\theta}$.
As seen from the Galois action on $\beta_{r,\theta}$, it is moreover
an $\mf{h}_r[G_{\Q}]$-module homomorphism
$\mcT_{r,\theta}/I_{r,\theta}\mcT_{r,\theta} \to (\mf{h}_{r,\theta}/I_{r,\theta})^{\,\sharp}(1)$.
The maps are compatible with $r$, and their projective limit is the desired surjective $\mf{h}[G_{\Q}]$-module homomorphism $\mcT_{\theta}/I_{\theta}\mcT_{\theta} \to Q_{\theta}$.

\end{sbpara}

\begin{sbpara} \label{unramified}

We explain how the surjection of \ref{Q} fits in an exact sequence
$$
	0\to P_{\theta} \to \mcT_{\theta}/I_{\theta}\mcT_{\theta} \to Q_{\theta} \to 0
$$
of $\mf{h}[G_{\Q}]$-modules that is canonically locally split over $G_{\qp}$.  

We use the fact that the Eisenstein part $\mcT^+_{\Eis} \to \mcT_{\quo,\Eis}$ of the canonical map of \ref{two_primes} is an isomorphism, or equivalently, that $\mcT_{\sub,\Eis} \to \mcT_{\Eis}/\mcT^+_{\Eis}$ is an isomorphism. To see this, one uses an $\mf{h}$-module splitting of the local exact sequence for $\mcT_{\theta}$ (see \cite{ohta-ord2}) and the method of Kurihara and Harder-Pink \cite{kurihara, hp}.  We refer the reader to \cite[Section 6.3]{fk-proof} for the argument.

Let us explain the use of this fact: by definition, complex conjugation acts on $Q_{\theta}$ by multiplication by $-1$.  Thus, $Q_{\theta}$ is a quotient of $\mcT_{\theta}/\mcT_{\theta}^+$.  By our isomorphism on Eisenstein components, it is a quotient of $\mcT_{\sub,\theta}/I_{\theta}\mcT_{\sub,\theta}$, which by \ref{ohta} is isomorphic to $\mf{h}_{\theta}/I_{\theta}$ as an $\mf{h}$-module.  This forces the quotient map to be an injection, so we have $Q_{\theta} \cong \mcT_{\sub,\theta}/I_{\theta}\mcT_{\sub,\theta}$.  But now, this tells us that $Q_{\theta}$ is an $\mf{h}[G_{\qp}]$-submodule of $\mcT_{\theta}/I_{\theta}\mcT_{\theta}$.  In other words, the surjection $\mcT_{\theta}/I_{\theta}\mcT_{\theta} \to Q_{\theta}$ is canonically locally split on $G_{\qp}$.  We then have necessarily that the kernel of the latter surjection is $\mcT_{\quo,\theta}/I_{\theta}\mcT_{\quo,\theta} \cong P_{\theta}$.   This yields the exact sequence.

It is perhaps worth observing that this sequence is also identified with the reduction modulo $I_{\theta}$ of the exact sequence of $\mf{h}[G_{\R}]$-modules
$0 \to \mcT^+_{\theta} \to \mcT_{\theta} \to \mcT_{\theta}/\mcT^+_{\theta} \to 0$.
Finally, the determinant of the $G_{\Q}$-action on $\mcT_{\theta}$ is known (e.g., from the determinants of modular Galois representations) and agrees with the $G_{\Q}$-action on $Q_{\theta}$, so the $G_{\Q}$-action on $P_{\theta}$ is trivial.

\end{sbpara}

\begin{sbpara}\label{upsilon}

We have that $P_{\theta}$ and $Q_{\theta}$ have trivial actions
of $G_{\tK}$.  Hence, we have a homomorphism
$$
	G_{\tK} \to \Hom_{\mf{h}}(Q_{\theta}, P_{\theta}), \qquad \sigma\mapsto (x\mapsto \sigma \tilde x - \tilde x),
$$
where $\tilde x$ is a lifting of $x$ to $\mcT_{\theta}/I_{\theta}\mcT_{\theta}$. 
By \ref{unramified}, this homomorphism factors through the unramified quotient $X$ of $G_{\tK}$. 
Thus, we have a homomorphism 
$X\to \text{Hom}_{\mf{h}}(Q_{\theta}, P_{\theta})$
that is compatible with the action of $\Gal(\tK/\Q)$. 
This gives a homomorphism of $\Gal(K/\Q)$-modules 
$$
	X^{-}(1) \to \text{Hom}_{\mf{h}}(Q_{\theta}(-1), P_{\theta}) \cong 
	\text{Hom}_{\mf{h}}((\mf{h}_{\theta}/I_{\theta})^{\,\sharp}, P_{\theta})\cong P_{\theta}^{\;\flat},
$$
where $P_{\theta}^{\;\flat}$ is $P_{\theta}$ on which $\sigma_a \in \Gal(\tK/\Q)$ 
acts as multiplication by $\langle a\rangle^{-1}$. In other words, we have a $\La_{\theta}$-module homomorphism
$$
	\Upsilon \colon Y_{\theta}=X^-(1)_{\theta} \to P_{\theta},
$$
with the Galois action of $G$ on the left and inverse diamond action of $G$ on the right.

\end{sbpara}

\begin{sbpara}\label{MW1} 

We describe the heart of the Mazur-Wiles proof of the Iwasawa main conjecture.

The Iwasawa main conjecture is the equality of ideals 
$$
	\cha_{\La_{\theta}}(Y_{\theta})=(\xi_{\theta}).
$$ 
By the analytic class number formula, this conjecture is reduced to $\cha_{\La_{\theta}}(Y_{\theta}) \subseteq (\xi_{\theta})$. 

Let $\mcL$ be the $\mf{h}[G_\Q]$-submodule of $\mcT_{\theta}$ generated by $\mcT_{\sub,\theta}$.  It follows as
in \ref{unramified} that we have an equality $\mcL_{\mf{m}} =\mcT_{\sub,\mf{m}}\oplus \mcL^+_{\mf{m}}$ of Eisenstein components.  Moreover, $P'_{\theta}=\mcL^+/I_{\theta}\mcL^+$ is $G_{\Q}$-stable in $\mcL/I_{\theta}\mcL$.  In other words, we have an exact sequence of $\mf{h}[G_{\Q}]$-modules
$$
	0\to P'_{\theta}\to \mcL/I_{\theta}\mcL\to Q_{\theta} \to 0.
$$ 
In the same way as $\Upsilon$, we may define $\Upsilon' \colon Y_{\theta} \to P'_{\theta}$, which is now surjective by construction.

The Iwasawa main conjecture can be deduced from this surjectivity of $\Upsilon'$. More precisely, we use the following facts:

\begin{enumerate}
	\item[(1)] The map $\Upsilon' \colon Y_{\theta} \to P'_{\theta}$ is surjective.
	\item[(2)] We have that $P'_{\theta}=\mcL^+/I_{\theta}\mcL^+$ with $\mcL^+$ a 
	finitely generated, faithful 
	$\mf{h}_{\theta}$-module.
	\item[(3)] The kernel of the canonical surjection $\La_{\theta}\to \mf{h}_{\theta}/I_{\theta}$ is contained in $(\xi_{\theta})$.\footnote{Actually, we know that the kernel coincides with $(\xi_{\theta})$, but this weaker statement is enough.}
\end{enumerate}
From (1), we obtain 
$$
	\cha_{\La_{\theta}}(Y_{\theta}) \subseteq \cha_{\La_{\theta}}(P'_{\theta}).
$$ 
From (2) and (3), we can deduce that
$$
	\cha_{\La_{\theta}}(P'_{\theta})\subseteq \cha_{\La_{\theta}}(\La_{\theta}/(\xi_{\theta}))=(\xi_{\theta}).
$$ 
Hence $\cha_{\La_{\theta}}(Y_{\theta})\subseteq (\xi_{\theta}).$

\end{sbpara}

\begin{sbpara}

The conjecture stated in the following section implies that $\Upsilon$ is surjective.  This tells us that the inexplicit lattice $\mcL$ required for the Mazur-Wiles proof in \ref{MW1} is precisely the canonical lattice $\mcT_{\theta}$.  In this sense, it suggests a refinement of the method of Ribet and Mazur-Wiles.   

\end{sbpara}

\subsection{The conjecture: $\varpi$ and $\Upsilon$ are inverse maps} \label{conj}

We state the conjecture of the third author \cite{sh-Lfn} and the result of the first two authors \cite{fk-proof}.

\begin{sbpara}

In \ref{varpi} and \ref{upsilon}, we  defined $\Lambda$-module homomorphisms
\begin{eqnarray*}
	\varpi \colon \mcP_{\theta} \to X^-(1)_{\theta} &\mr{and}& \Upsilon \colon X^-(1)_{\theta} \to \mcP_{\theta}.
\end{eqnarray*}
We have the conjecture of the third author.  See \cite[Conjecture 4.12]{sh-Lfn}, where the conjecture is given up to a canonical unit; this stronger version was a stated hope of the third author.

\begin{conj}
	The maps $\varpi$ and $\Upsilon$ are inverse to each other.
\end{conj}

This conjecture provides an explicit description of $X^-(1)_{\theta}$ in terms of modular symbols.  In this sense, it may be viewed as a refinement of the main conjecture.

\end{sbpara}

\begin{sbpara} \label{MW2}  

We state the result \cite[Theorem 7.2.3(1)]{fk-proof} of the first two authors.  Let $\xi'_{\theta} \in \La_{\theta}$ denote the derivative of the $p$-adic $L$-function $\xi_{\theta}$ in the $s$-variable (see \ref{xi}).

\begin{thm}
	We have $\xi'_{\theta}\Upsilon \circ \varpi = \xi'_{\theta}$ modulo $p$-torsion in $\mcP_{\theta}$. 
\end{thm}

If $\xi_{\theta}$ has no multiple roots, the theorem implies the conjecture up to $p$-torsion in $P_{\theta}$.  
In fact, it leads to proofs of the conjecture under various hypotheses: see \cite[Section 7.2]{fk-proof}.  

\end{sbpara}

\begin{sbpara}

McCallum and the third author conjectured that the image of the cup product
$$
	H^1(\Z[\zeta_{Np^r},\tfrac{1}{Np}],\zp(1)) \otimes_{\zp} H^1(\Z[\zeta_{Np^r},\tfrac{1}{Np}],\zp(1))
	\xrightarrow{\cup} H^2_{\et}(\Z[\zeta_{Np^r},\tfrac{1}{Np}],\zp(2))
$$	
projects onto $H^2_{\et}(\Z[\zeta_{Np^r},\tfrac{1}{Np}],\zp(2))_{\theta}^+$ 
(\cite{mcs} for $N = 1$), which implies that $\varpi$ is surjective. This generation conjecture follows if we know that $\varpi \circ \Upsilon = 1$.  In particular, it holds if $\xi_{\theta}$ has no multiple roots, and 
it also holds if $\mcP_{\theta} \otimes_{\zp} \qp$ is generated by one element over $\Lambda_{\theta} \otimes_{\zp} \qp$ 
\cite[Theorem 7.2.8]{fk-proof}.   

\end{sbpara}

\subsection{The proof that $\xi'\Upsilon \circ \varpi = \xi'$}\label{key_ideas}

We explain some of the important aspects of the proof of the main theorem, referring to the relevant sections of \cite{fk-proof} for details.

\begin{sbpara}\label{newdiag}  We consider a refinement of the diagram in 
\ref{key_diag} in which we
divide the right-hand square of that diagram into two squares:
$$
	\SelectTips{eu}{} 
	\xymatrix{ \mcS_{\theta} \ar[r]^-z \ar[d]^{\bmod I} & \vpr H^2_{\et}(Y_1(Np^r), \Z_p(2))_{\theta} 
	\ar@{-->}[r]^-{\HS} \ar[d]^{\infty}
	& H^1_{\et}(\Z[\tfrac{1}{Np}], \mcT_{\theta}(1)) \ar[r]^-{\text{reg}} \ar[d] & \mf{S}_{\theta} 
	\ar[d]^{\bmod I} \\
	P_{\theta} \ar[r]^{\varpi} & Y_{\theta} \ar[r]^{\xi'_{\theta}}  & Y_{\theta} \ar[r]^{\Upsilon} & P_{\theta}.} 
$$
Here, the maps $z$ and $\infty$ are the projective limits of the $\theta$-components of the maps $z_r$ and $\infty_r$.  The commutativity of the left square of the diagram in \ref{newdiag} is seen in Section \ref{zeta_elts}.   The discussion of the rest of this diagram, and the fact that the bottom row is also multiplication by $\xi'_{\theta}$, compose the rest of this subsection.  

It is remarkable that $\xi'_{\theta}$ appears here in two very different contexts.  The $\xi'_{\theta}$ that appears in the diagram and contributes to $\xi'_{\theta} \Upsilon \circ \varpi$ is related to cup product with the logarithm of the cyclotomic character.  The other $\xi'_{\theta}$ is the constant term modulo $\xi_{\theta}$ of a $\La$-adic modular form that appears in a computation of the regulators of zeta elements.

\end{sbpara}

\begin{sbpara} \label{HS}
	The map $\HS$ arises from the Hochschild-Serre spectral sequences
	$$
		E_2^{i,j} = H^i(\Z[\tfrac{1}{Np}],H^j(Y_1(Np^r)_{/\qbar},\zp(2))) \Rightarrow 
		E^{i+j} = H^{i+j}(Y_1(Np^r),\zp(2)).
	$$
	as the projective limit over $r$ of maps $E^2 \to E_2^{1,1}$, followed by projection to the ordinary
	$\theta$-part.  We remark that
	$$
		H^1_{\et}(\Z[\tfrac{1}{Np}], \mcT_{\theta}(1)) \subset 
		H^1_{\et}(\Z[\tfrac{1}{Np}], \tT_{\theta}(1)),
	$$
	and the image of $\HS$ is actually contained in the larger group, hence the dotted arrow.  
	However, elements of $\mcS_{\theta}$
	are carried to the smaller group under $\HS \circ z$ \cite[Proposition 3.3.14]{fk-proof}, 
	so we can still make sense of the diagram.
\end{sbpara}
	
\begin{sbpara}\label{third_arrow}

	The third vertical arrow in the diagram of \ref{newdiag} is a composition of maps as follows:
	$$
		H^1_{\et}(\Z[\tfrac{1}{Np}], \mcT_{\theta}(1))\to 
		H^1_{\et}(\Z[\tfrac{1}{Np}], Q_{\theta}(1))\xrightarrow{\cup\, (1-p^{-1})\log(\kappa)} 
		H^2_{\et}(\Z[\tfrac{1}{Np}], Q_{\theta}(1)) \xsim Y_{\theta}.
	$$	
	The first map is induced by the surjection $\mcT_{\theta} \to Q_{\theta}$.  
	The map $\kappa \colon G_{\Q} \to \zp^{\times}$ is the $p$-adic cyclotomic character,
	and $\log$ is the $p$-adic logarithm.  The second map is the cup product, where we regard
	$\log(\kappa) = \log \circ \kappa$ as an element of $H^1_{\et}(\Z[\tfrac{1}{p}], \Z_p)$.
	
	For the third map, note that $Q_{\theta} \cong (\Lambda_{\theta}/\xi_{\theta})^{\,\sharp}(1)$
	by \ref{hmodI} and \ref{Q}.
	As the $p$-cohomological dimension of $\Z[\frac{1}{Np}]$ is $2$, the group 
	$H^2_{\et}(\Z[\frac{1}{Np}],(\Lambda_{\theta}/\xi_{\theta})^{\,\sharp}(2))$ is isomorphic to the 
	quotient of 
	$$
		H^2_{\et}(\Z[\tfrac{1}{Np}],\Lambda_{\theta}^{\,\sharp}(2))_{\theta} \xsim \,
		\vpr H^2_{\et}(\mc{O}_r[\tfrac{1}{Np}],
		\zp(2))_{\theta} \xsim Y_{\theta}
	$$
	by the $\Gal(K/\Q)$-action of $\xi_{\theta}$.  Here, the first isomorphism is by Shapiro's lemma, and
	the second is from \ref{varpi}.  
	By the main conjecture 
	and the fact that $Y_{\theta}$ has no finite $\La$-submodules, $Y_{\theta}$ is $\xi_{\theta}$-torsion,
	so $Y_{\theta}/\xi_{\theta}Y_{\theta} = Y_{\theta}$.  
	Putting this all together, we have the map.
  
\end{sbpara}

\begin{sbpara}\label{D}

We define a functor $D$ on pro-$p$ $G_{\qp}$-modules.

Let $T = \varprojlim_{\lambda} T_{\lambda}$ for a projective system of finite abelian $p$-groups 
$T_{\lambda}$.  For any abelian group $M$, set 
$T \cotimes M = \varprojlim_{\lambda}(T_{\lambda} \otimes M)$.
Let $W$ denote the Witt vectors of $\overline{\F}_p$, and suppose that
the $T_{\lambda}$ are endowed with compatible actions of $h[G_{\qp}]$ for a pro-$p$ ring $h$.
We may then consider the $h$-module $D(T)$ that is the fixed part 
$$
	D(T)  = (T \cotimes W)^{G_{\qp}}
$$ 
for the diagonal action of $G_{\qp}$ on $T \cotimes W$.  If the $G_{\qp}$-actions on the $T_{\lambda}$ are unamified, then $D(T)$ and $T$ are isomorphic $h$-modules.  If $T$ has trivial $G_{\qp}$-action,
then $D(T) \cong T \cotimes W^{G_{\qp}} = T \cotimes \zp \cong T$, and this isomorphism is canonical.
See \cite[Section 1.7]{fk-proof}.

\end{sbpara}

\begin{sbpara} \label{preg}

We define $p$-adic regulator maps for unramified, pro-$p$ $G_{\qp}$-modules.

Let $T$ be as in \ref{D}, and suppose that the action of $G_{\qp}$ on $T$ is unramified.
Let $E = Q(W)$ be the maximal unramified extension of $\qp$. 
The $p$-adic regulator map \cite[Section 4.2]{fk-proof}
$$
	\reg_T \colon H^1_{\et}(\Q_p, T(1))\to D(T)
$$
for $T$ is the $h$-module homomorphism defined as the composition
$$
	H^1_{\et}(\Q_p, T(1)) \xrightarrow{\mr{inf}} H^1_{\et}(E, T(1))^{\text{Fr}_p=1} \xsim 
	(T \cotimes E^\times)^{\Fr_p=1} \to D(T).
$$
Here, the first map is inflation, the second is Kummer theory, and the final map is induced by 
$$
	E^{\times} \to W(\overline{\F}_p), \quad x \mapsto p^{-1}\log\Big(\frac{x^p}{\Fr_p(x)}\Big),
$$
where the $p$-adic logarithm $\log$ is defined to take $p$ to $0$.  

Note that if $G_{\qp}$ acts trivially on $T$, then $\reg_T$ is induced by the map $(1-p^{-1})\log \colon \qp^{\times} \to \zp$ in a similar fashion.

\end{sbpara}

\begin{sbpara}\label{ES}

We define the $p$-adic regulator map $\reg$ in the diagram of \ref{newdiag}.

Note that $\mcT_{\quo}$ has by definition an unramified $G_{\qp}$-action.
We have a refinement \cite[Section 1.7]{fk-proof} of Ohta's $\La$-adic Eichler-Shimura isomorphism \cite{ohta-eich}.  That is, there is a canonical isomorphism of $\mf{h}$-modules $D(\mcT_{\quo}) \xsim \mf{S}$, and in particular $\mcT_{\quo}$ and $\mf{S}$ are noncanonically isomorphic.  The map $\reg$ is then defined as the composition 
$$
	\reg \colon H^1_{\et}(\Z[\tfrac{1}{Np}], \mcT_{\theta}(1))\to H^1_{\et}(\Q_p, \mcT_{\quo,\theta}(1)) 
	\xrightarrow{\reg_{\mcT_{\quo}}} D(\mc{T}_{\quo}) \xsim \mf{S}_{\theta}.
$$ 

\end{sbpara}

\begin{sbpara}  \label{3S_isos}

We explain the right-hand vertical map ``mod $I$'' in the diagram of \ref{newdiag}.

The $G_{\qp}$-action on $\mcT_{\quo}/I\mcT_{\quo}$ is trivial, so the canonical isomorphism $D(\mcT_{\quo}) \xsim \mf{S}$ provides an isomorphism 
$\mcT_{\quo}/I\mcT_{\quo} \xsim \mf{S}/I\mf{S}$.  In particular, we obtain ``mod $I$'' as the 
composition of projection followed by a string of canonical isomorphisms:
$$
	\mf{S}_{\theta} \twoheadrightarrow 
	\mf{S}_{\theta}/I_{\theta}\mf{S}_{\theta} \xsim \mcT_{\quo,\theta}/I_{\theta}\mcT_{\quo,\theta}
	\xsim \mcT^+_{\theta}/I_{\theta}\mcT^+_{\theta} \xsim
	\mcS_{\theta}/I_{\theta}\mcS_{\theta} = P_{\theta}. 
$$

\end{sbpara}

\begin{sbpara} 

The commutativity of the two right-hand squares in the diagram of \ref{newdiag} are nontrivial cohomological exercises.  We mention only which calculations must in the end be performed.

\begin{enumerate}
	\item[(1)]
	The commutativity of the middle square is reduced to that (see \cite[Section 9.4]{fk-proof}) of
	$$
		\SelectTips{eu}{} \xymatrix@C=60pt{
		H^1_{\et}(\Z[\tfrac{1}{Np}], (\Lambda_{\theta}/\xi_{\theta})^{\,\sharp}(2)) 
		\ar[r]^{\cup\, (1-p^{-1})\log(\kappa)} \ar[d]^{\wr} & 
		H^2_{\et}(\Z[\tfrac{1}{Np}], (\Lambda_{\theta}/\xi_{\theta})^{\,\sharp}(2)) \\
		H^2_{\et}(\Z[\tfrac{1}{Np}], \Lambda^{\,\sharp}_{\theta}(2)) \ar[r]^{\xi'_{\theta}} & 
		H^2_{\et}(\Z[\tfrac{1}{Np}], \Lambda^{\,\sharp}_{\theta}(2)), \ar[u]^{\wr} }
	$$
	the vertical arrows occurring in the long exact sequence in the $\Z[\frac{1}{Np}]$-cohomology of 
	$$
		0 \to \Lambda^{\,\sharp}_{\theta}(2) \to \Lambda^{\,\sharp}_{\theta}(2) \to 
		(\Lambda_{\theta}/\xi_{\theta})^{\,\sharp}(2) \to 0.
	$$
	Thus, the $\xi'_{\theta}$ that appears in the diagram is found in Galois cohomology.
	\item[(2)]
	The commutativity of the right-hand square is reduced to verifying that
	the map 
	$$
		Y_{\theta} \cong H^2_{\et}(\Z[\tfrac{1}{Np}],Q_{\theta}(1)) \leftarrow 
		H^2_{\et}(\Z[\tfrac{1}{Np}],\mc{T}_{\theta}/I_{\theta}\mc{T}_{\theta}(1)) \to H^2_{\et}(\qp,P_{\theta}(1))
		\cong P_{\theta}
	$$
	given by lifting and then projecting is well-defined and agrees with $\Upsilon$ \cite[Section 9.5]{fk-proof}.  
	Here, the first isomorphism was discussed in \ref{third_arrow} and the last is the invariant map of 
	local class field theory, recalling from \ref{unramified} that $P_{\theta}$ has trivial Galois action.  This
	description is closer to the construction of $\Upsilon$ that will appear for $\F_q(t)$ in Section \ref{F_q(t)}.
\end{enumerate}

\end{sbpara}

\begin{sbpara} 

It remains to prove that the composition $\mcS_{\theta}/I_{\theta}\mcS_{\theta} \to \mf{S}_{\theta}/I_{\theta}\mf{S}_{\theta} \to P_{\theta}$, where the first arrow is the composition of the upper horizontal arrows modulo $I_{\theta}$ in the diagram of \ref{newdiag}, coincides with multiplication by $\xi'_{\theta}$ on $P_{\theta}$.  This is deduced in \cite[Sections 4.3 and 8.1]{fk-proof} from the computation of the $p$-adic regulators of zeta elements given in \cite{Oc, Fu}.  This is a very delicate analysis: we explain only the rough idea of how $\xi'_{\theta}$ appears at its end.  

The map $P_{\theta} \to P_{\theta}$ is shown to be given (modulo $\xi_{\theta}$) by multiplication by the constant term at $t = 1$ of a $p$-adic $L$-function in a variable $t$ that takes values in $\mf{M}_{\theta}$.  This $p$-adic $L$-function is a product of two $\Lambda$-adic Eisenstein series which vary with $t$.  The constant term in the $q$-expansion of this product is itself a product of two zeta functions $\zeta_p(t)\xi_{\theta}(s+t-1)$, where $\zeta_p(t)$ is the $p$-adic Riemann zeta function and $s$ is the variable for $\Lambda_{\theta} \subset \mf{h}_{\theta}$.  Note that $\zeta_p(t)$ has a simple pole at $t = 1$ with residue $1$.  To evaluate $\zeta_p(t)\xi_{\theta}(s+t-1)$ modulo $\xi_{\theta}(s)$ at $t = 1$, we can first subtract $\zeta_p(t)\xi_{\theta}(s)$ from the product and then take the resulting limit 
$$
	\lim_{t \to 1} \frac{\xi_{\theta}(s+t-1) - \xi_{\theta}(s)}{t-1} = \xi'_{\theta}(s).
$$
In this manner, the map is shown to be multiplication by $\xi'_{\theta}$.

\end{sbpara}

\section{The case of $\GL_2$ over $\F_q(t)$} \label{F_q(t)}

We now consider the field $F = \F_q(t)$ for some prime power $q$.  In this section, we provide $F$-analogues of the constructions, conjecture, and theorem of Section
\ref{GL_2-Q}.  We require the following objects:
\begin{enumerate}
	\item[$\bullet$] the ring $\mc{O} = \F_q[t]$,  
	\item[$\bullet$] the completion $F_{\infty}=F((t^{-1}))$ of $F$ at the place $\infty$,
	\item[$\bullet$] the valuation ring $\mc{O}_{\infty} = \F_q\ps{t^{-1}}$ of $F_{\infty}$, which does
	not contain $\mc{O}$,
	\item[$\bullet$] a prime number $p$ different from the characteristic of $\F_q$,
	\item[$\bullet$] a non-constant polynomial $N \in \mc{O}$.  
\end{enumerate}	
Let us also fix an embedding $\Fbar \hookrightarrow \overline{F_{\infty}}$ of separable closures.
To avoid technical complications, we assume in this section that $p$ does not divide 
 $(q+1)|(\mc{O}/N\mc{O})^{\times}|$.

The organization of this section follows closely that of Section \ref{GL_2-Q}.
We hope to make clear that most constructions are remarkably similar to the case of $\Q$, though we also highlight differences.  We work with congruence subgroups of $\GL_2(\mc{O})$, rather than of $\SL_2(\Z)$.
Modular symbols, used to construct $\varpi$, are now found in the homology $\mcS$ of the compactification of the quotient of the Bruhat-Tits tree by a congruence subgroup.  This $\mcS$ is a quotient of an \'etale homology group $\mcT$ of the Drinfeld modular curve, used in constructing $\Upsilon$.  As most constructions are so similar, we provide less detail than in Section \ref{GL_2-Q}.  We intend for full details to appear in a forthcoming paper.

\subsection{From modular symbols to cup products: the map $\varpi$} 

\begin{sbpara} \label{homol_F}

We introduce homology groups $\mcS$ and $\tS$ of the Bruhat-Tits tree.
	
Consider the Bruhat-Tits tree $B$ for $\PGL_2(F_{\infty})$.  Its vertices are homothety classes of $\mc{O}_{\infty}$-lattices $\mcL$ of rank $2$ in $F_{\infty}^2$, or equivalently, elements of $\PGL_2(F_{\infty})/\PGL_2(\mc{O}_{\infty})$.  This tree is $(q+1)$-valent, and two lattices $\mcL \subset \mcL'$ connected by an edge if $[\mcL' : \mcL] = q$.\footnote{Note that $q$ appears in this sentence as the order of the residue field of $\mc{O}_{\infty}$.}  
The oriented edges then correspond to elements of $\PGL_2(F_{\infty})/\mc{I}_{\infty}$, where
$\mc{I}_{\infty}$ is the Iwahori subgroup of matrices in $\PGL_2(\mc{O}_{\infty})$ that are upper-triangular modulo the maximal ideal of $\mc{O}_{\infty}$.
The group $\PGL_2(F)$  acts on the left on $B$ in the evident manner.

Let $\tGa_1(N)$ be the congruence subgroup of $\GL_2(\mc{O})$ given by
$$
	\tGa_1(N) = \left\{\begin{pmatrix} a & b \\ c & d\end{pmatrix} \in \GL_2(\mc{O}) \;\Big|\; (c,d) \equiv
	(0,1) \bmod N \right\}.
$$
We may complete the Bruhat-Tits tree to a space $B^*$ by adding in the (rational) ends, which correspond to elements of $\mb{P}^1(F)$.
We define
\begin{eqnarray*}
	U(N) = \tGa_1(N) \bs B &\mr{and}& \overline{U}(N)  = \tGa_1(N) \bs B^*.
\end{eqnarray*}
The elements of $\tGa_1(N)\bs \mb{P}^1(F)$ are the ends of $U(N)$. 
Our homology groups, or spaces of modular symbols, are then
$$
	\mcS = H_1(\overline{U}(N),\zp) \subset \tS =H_1(\overline{U}(N),\{\mr{ends}\},\zp).
$$

\end{sbpara}

\begin{sbpara}\label{3.2.3}  

We introduce Manin-Teitelbaum symbols $[u:v] \in \tS$.

Modular symbols in $\tS$  were defined by Teitelbaum \cite{teitelbaum} analogously to the case of $\Q$.  In particular, given $\alpha, \beta \in \mb{P}^1(F)$, we have a modular symbol that is the class $\{\alpha \to \beta\}$ of any non-backtracking path in the Bruhat-Tits tree that connects the two corresponding ends of $B$.

Analogues of Manin symbols are defined as before.  That is, for $u, v \in \mc{O}/N\mc{O}$ with $(u,v) = (1)$, we choose $\gamma = \left(\begin{smallmatrix} a & b \\ c & d\end{smallmatrix}\right) \in \GL_2(\mc{O})$ with $u = c \bmod N$ and $v = d \bmod N$, and then 
$$
	[u:v] = \left\{ \frac{d}{bN} \to \frac{c}{aN}\right\}.
$$  
These symbols generate $\widetilde{H}$ and yield a presentation with identical relations to those of \ref{maninsym}.

\end{sbpara}

\begin{sbpara}

We introduce the intermediate space $\tS^0$ on which we define $\varpi$.  

Let $\tS^0$ denote the $\zp$-submodule of $\tS$ generated by the Manin symbols $[u:v]$ with $u,v \neq 0$.  As in the case of $\GL_2$ over $\Q$, we have $\mcS \subset \tS^0 \subset \tS$.

\end{sbpara}

\begin{sbpara} \label{cycl_units}

We introduce cyclotomic $N$-units $\la{u}$ in abelian extensions $F_N \subset E_N$ of $\F_q(t)$.  The reader may find a powerful analogy with objects in the theory of cyclotomic fields over $\Q$.

We consider the cyclotomic $N$-units $\la{u}$ for $u \in \mc{O} - (N)$. 
These are the roots of the Carlitz polynomials \cite{carlitz} for divisors of $N$, or are equivalently the $N$-torsion points of the Carlitz module.  As $\la{u}$ depends only on $u$ modulo $N$, we abuse notation and consider it for nonzero $u \in \mc{O}/N\mc{O}$. 
We can visualize $\la{u}$ in  the completion $\C_{\infty}$ of $\overline{F_{\infty}}$
by
$$
	\la{u} =  \exp\left(\frac{u\pi}{N}\right) = \frac{u}{N} \prod_{a \in \mc{O} - \{0\}} 
	\left( 1 - \frac{u}{Na} \right),
$$
where $\exp$ is the Carlitz exponential and $\pi \in \C_{\infty}$ is transcendental over $F$.  

Let $E_N = F(\la{1})$, which is an abelian extension of $F$ of conductor $N\infty$ containing
no constant field extension of $F$.  
There is an isomorphism $\Gal(E_N/F) \xsim (\mc{O}/N\mc{O})^{\times}$ such that
$a \in (\mc{O}/N\mc{O})^{\times}$ is the image of an element $\sigma_a \in \Gal(E_N/F)$ that satisfies $\sigma_a(\la{1}) = \la{a}$.
Let $F_N$ be the largest subfield of $E_N$ in which $\infty$ splits completely over $F$, which we might call 
the ray class field of modulus $N$.  Under the above isomorphism, $\Gal(E_N/F_N)$ is identified with $\F_q^{\times}$.  In fact, 
we have $\sigma_c(\la{u}) = c\la{u}$ for $c \in \F_q^{\times}$. These facts are found in the work of Hayes \cite{hayes}.

Let $\mc{O}_N$ denote the integral closure of $\mc{O}$ in $F_N$.
Since $p \nmid (q-1)$ by assumption, the image of $\la{u}$ in the $p$-completion of the $N$-units of $E_N$ is fixed by the action of $\F_q^{\times}$.  This allows us to view $\la{u}$ as an element of $H^1_{\et}(\mc{O}_N[\frac{1}{N}], \Z_p(1))$.  For nonzero $u,v\in \mc{O}/N\mc{O}$, we may consider the cup product 
$$
	\la{u}\cup  \lambda_{\frac{v}{N}} \in H^2_{\et}(\mc{O}_N[\tfrac{1}{N}], \Z_p(2)).
$$

\end{sbpara}

\begin{sbpara}  We define the map $\varpi$.  Here, we work directly with \'etale cohomology, rather than $K_2$.

There is a homomorphism 
$$
	\varpi \colon \tS^0 \to H^2_{\et}(\mc{O}_N[\tfrac{1}{N}], \Z_p(2)), \quad 
	[u:v] \mapsto \la{u} \cup \lambda_{\frac{v}{N}}.
$$
In the current setting, we can no longer quickly verify from the presentation of $\mc{M}^0$ 
that $\varpi$ is well-defined.  Rather, we see this as a consequence of the
argument that $\varpi$ is ``Eisenstein'' in Section \ref{zeta_F_q(t)}.
 
\end{sbpara} 
 
\begin{sbpara}
We introduce the cuspidal Hecke algebra $\mf{h}$ and its Eisenstein ideal $I$.

Let $\mf{n}$ denote a nonzero ideal of $\mc{O}$.
Through the action of $\GL_2(F)$ on $B$, we have a Hecke operator $T(\mf{n})$ acting on $\mcS$
as the correspondence associated to $\tGa_1(N) \left( \begin{smallmatrix} 1&0\\0&n \end{smallmatrix} \right)\tGa_1(N)$, where $\mf{n} = (n)$.  Let $\mf{h}$ be the subring of $\End_{\Z_p}(\mcS)$ generated over $\Z_p$ by the Hecke operators $T(\mf{n})$.  

We also have diamond operators $\langle \mf{a} \rangle$ in $\mf{h}$ for nonzero ideals $\mf{a}$ of $\mc{O}$ prime to $(N)$.  This $\langle \mf{a} \rangle$ depends only on the reduction modulo $N$ of the monic generator of $\mf{a}$.  

The Eisenstein ideal $I$ is the ideal of $\mf{h}$ generated by $T(\mf{n}) - \sum_{\mf{d} \mid \mf{n}} \mf{N}(\mf{d})\langle \mf{d}\rangle$ for all nonzero ideals $\mf{n}$ of $\mc{O}$, taking $\langle \mf{d} \rangle = 0$ if $\mf{d}+(N) \neq (1)$.  Here, $\mf{N}(\mf{n}) = [\mc{O}:\mf{n}]$ is the absolute norm of $\mf{n}$.  

Similarly, we have the Eisenstein ideal $\mf{I}$ of the Hecke algebra $\mf{H} \subset \End_{\zp}(\tS)$.
  
\end{sbpara} 

\begin{sbpara}

To say that $\varpi$ is ``Eisenstein'' is to say that $\varpi$ factors through a map
$$
	\varpi \colon \tS^0/\mf{I}\tS^0 \to H^2_{\et}(\mc{O}_N[\tfrac{1}{N}],\zp(2)).
$$
We explain this result in Section \ref{zeta_F_q(t)}. 

\end{sbpara}
 
\begin{sbpara}\label{3.1.11} Let $G=(\mc{O}/N\mc{O})^\times/\F_q^\times$, and set $\Lambda = \zp[G]$. 

\begin{enumerate}
	\item[(1)] We have a ring homomorphism $\iota \colon \Lambda \to \mf{h}$ which sends the group 
	element $[a]$ in $\Z_p[G]$ for $a \in G$ to the inverse $\langle a \rangle^{-1}$ of the diamond operator 
	corresponding to $a$.
	\item[(2)] We have the isomorphism $\Gal(F_N/F) \xsim G$ of class field theory (see
	\ref{cycl_units}).
\end{enumerate}
Modules over $\mf{h}$ and $\zp[\Gal(F_N/F)]$
become $\Lambda$-modules through these identifications.

\end{sbpara}

\subsection{Working with fixed level} \label{notower}

We explain why we work with fixed level in Section \ref{F_q(t)}, and we define our two objects of study.

\begin{sbpara}

We do not pass up a tower for the following reason on the $\GL_1$-side.  By assumption on $p$, the field $\F_q$ has no nontrivial $p$th roots of unity.  Since $F_N/F$ contains no constant field extension, $F_N$ also contains no nontrivial $p$th roots of unity.  So, even if we ``increase'' $N$, we are unable to employ the Iwasawa-theoretic trick of passing Tate twists through projective limits of Galois cohomology groups.  In particular, since we deal with cohomology with $\zp(2)$-coefficients, we do not work with class groups.

\end{sbpara}

\begin{sbpara} \label{two_obj_F_q(t)} 

We again have two objects of study:
\begin{enumerate}
	\item[$\bullet$] the geometric object $P = \mc{S}/I\mc{S}$ for $\GL_2$,
	\item[$\bullet$] the arithmetic object $Y=H^2_{\et}(\mc{O}_N, \Z_p(2))$ for $\GL_1$. 
\end{enumerate}

Given a character $\theta \colon G \to \qpbt$, we set $\Lambda_{\theta} = \zp[\theta]$ and view it 
as a quotient of $\Lambda$ through $\theta$.   For a $\Lambda$-module $M$, we let $M_{\theta} = M \otimes_{\Lambda} \Lambda_{\theta}$ denote the $\theta$-part of $M$.  If $\theta$ is primitive, then 
our assumption that $p$ does not divide $|G|$ implies that 
the canonical maps
\begin{eqnarray*}
	P_{\theta} \to \tS^0_{\theta}/\mf{I}_{\theta}\tS^0_{\theta} &\mr{and}& Y_{\theta} \to 
	H^2_{\et}(\mc{O}_N[\tfrac{1}{N}], \Z_p(2))_{\theta}
\end{eqnarray*}
are isomorphisms.

\end{sbpara}

\subsection{Zeta elements: $\varpi$ is ``Eisenstein''} \label{zeta_F_q(t)}

We explain that $\varpi$ factors through the quotient of $\tS^0$ by the Eisenstein ideal $\mf{I}$.

\begin{sbpara}
	
We define Siegel units on Drinfeld modular curves.  

Let $Y(N)$ denote the Drinfeld modular curve that is the moduli scheme for pairs consisting of a rank $2$ Drinfeld module over an $\mc{O}[\tfrac{1}{N}]$-scheme and a full $N$-level structure (or, basis of the $N$-torsion) on it. 
Over $Y(N)$, we have a universal Drinfeld module, equipped with a full $N$-level structure, which locally looks like $(N^{-1}\mc{O}/\mc{O})^2 \times Y(N)$.  On the universal Drinfeld module is a certain theta function $\Theta$.  Given an element of $(\tfrac{u_1}{N},\tfrac{u_2}{N}) \in (N^{-1}\mc{O}/\mc{O})^2$, we may pull $\Theta$ back to a unit on the Drinfeld modular curve using the second coordinate of the level structure.  This unit $g_{\alpha,\beta} \in \mc{O}_{Y(N)}^{\times}$ is the analogue of a Siegel unit.\footnote{Actually, $g_{\alpha,\beta}$ as we have described it is not well-defined until we take its $q^2-1$ power.  The assumption that $p \nmid (q^2-1)$ is used to avoid this issue when we work with \'etale cohomology.} 

Let $Y_1(N)$ be the moduli scheme for pairs consisting of a rank $2$ Drinfeld module over an $\mc{O}[\frac{1}{N}]$-scheme and a point of order $N$ on it.  If we take $\alpha = 0$, then the Siegel unit $g_{0,\beta}$ may again be viewed as an element of $\mc{O}_{Y_1(N)}^{\times}$. 	

\end{sbpara}

\begin{sbpara}

If we take a $K$-theoretic product of two Siegel-type units, we obtain the Beilinson-type elements considered by Kondo and Yasuda \cite{ky}.  See also the work of Kondo \cite{kondo} and Pal \cite{pal}.  
Much as in the case of $\Q$, we have a map
$$
	z \colon \tS^0 \to H^2_{\et}(Y_1(N),\zp(2)), \qquad [u:v] \mapsto 
	g_{0,\tfrac{u}{N}} \cup g_{0,\tfrac{v}{N}}
$$
of $\mf{H}$-modules.
We can specialize this at the cusp corresponding to $\infty \in \mb{P}^1(F)$ to obtain $\la{u} \cup \lambda_{\frac{v}{N}}$.  This specialization map $\infty$ is Eisenstein.  Hence, we see that 
$$
	\varpi = \infty \circ z \colon \tS^0 \to H^1_{\et}(\mc{O}_N[\tfrac{1}{N}],\zp(2))
$$ 
is well-defined and Eisenstein.

\end{sbpara}

\subsection{Homology of Drinfeld modular curves}

In this subsection, we study the \'etale homology groups of Drinfeld modular curves.  Unlike in Section \ref{ordinary}, we do not take ordinary parts.  That is, the Galois representations found in the homology of Drinfeld modular curves are already ``special at $\infty$,'' the required analogue of ``ordinary at $p$.''   Moreover, the resulting unramfied-at-$\infty$ quotient may in the present setting be identified with the space $\mcS$ of cuspidal symbols, which is the analogue of the plus quotient of homology of \ref{limitpi}.   In other words, the place $\infty$ of $\F_q(t)$ plays both the roles that $p$ and the real place do in the $\GL_2$-setting over $\Q$. 

The statements in this subsection are consequences of the work of Drinfeld \cite{drinfeld-ell}.

\begin{sbpara}\label{3.1.1}
We first introduce the \'etale homology group $\mcT$.

Over $F$, the Drinfeld modular curve $Y_1(N)_{/F}$ has a smooth compactification
$X_1(N)_{/F}$.   Over $\C_{\infty}$ (or $\Fbar$), it is given by adding in the set of cusps
$\tGa_1(N) \backslash \mb{P}^1(F)$ of the Drinfeld upper half-plane.   We define our \'etale homology group
$$
	\mcT = H_1^{\et}(X_1(N)_{/\Fbar},\zp)
$$
as the $\Z_p$-dual of $H^1_{\et}(X_1(N)_{/\Fbar}, \Z_p)$.

The Hecke algebra generated by the $T(\mf{n})$ in $\End_{\zp}(\mcT)$ is in fact equal to $\mf{h}$.  The module 
$\Q_p\otimes_{\Z_p} \mcT$ over the total quotient ring $\qp \otimes_{\zp} \mf{h}$ of $\mf{h}$ is free of rank $2$.

 \end{sbpara}

\begin{sbpara} \label{one_prime}
 
We study the action of $G_{F_{\infty}}$ on $\mcT$.  

We have an exact sequence 
$$
	0 \to \mcT_{\sub} \to \mcT \to \mcT_{\quo} \to 0
$$
of $\mf{h}[G_{F_{\infty}}]$-modules, with $\mcT_{\sub}$ and $\mcT_{\quo}$ defined as follows.  First, 
$\mcT_{\sub}$ is the largest submodule of $\mcT$ such that $G_{F_{\infty}}$ acts on $\mcT_{\sub}(-1)$
trivially, and $\mcT_{\quo}$ is the quotient.  Then $\mcT_{\quo}$ is equal to the maximal unramified, 
$\mf{h}$-torsion-free quotient of $\mcT$.  In this way, the place $\infty$ plays the role that the place at $p$ does in 
\ref{two_primes}.  In fact, $G_{F_{\infty}}$ acts trivially on $\mcT_{\quo}$, and both $\qp \otimes_{\zp} \mcT_{\sub}$ and $\qp \otimes_{\zp} \mcT_{\quo}$ are free of rank $1$ over $\qp \otimes_{\zp} \mf{h}$.  

The above short exact sequence is split as a sequence of $\mf{h}$-modules: $\mcT_{\quo}$ is the isomorphic image of the $\mf{h}$-submodule of $\mcT$ on which a choice of Frobenius element acts trivially.  This will be used in constructing $\Upsilon$ below.

\end{sbpara}

\begin{sbpara} \label{TU} 

Since $\overline{U}(N)$ is essentially the graph of the special fiber of a model of $X_1(N)$ over $\mc{O}_{\infty}$,
we have a surjective homomorphism 
$$\mcT\to \mcS = H_1(\overline{U}(N), \Z_p).$$
Via this map, $\mcS$ is identified with the quotient $\mcT_{\quo}$ of $\mcT$ with trivial $G_{F_{\infty}}$-action.  
In this way, $\mcT_{\quo}$ is also
analogous to the plus quotient of homology in \ref{limitpi}.  That is, the place $\infty$ also plays the role 
that the real place does over $\Q$.

\end{sbpara}

\begin{sbpara}

Let $\mf{S}$ be the space of those $\zp$-valued, special-at-$\infty$ cuspidal automorphic forms 
$$
	\phi \colon \PGL_2(F) \backslash \PGL_2(\mb{A}_F) / (K_f \times \mc{I}_{\infty}) \to \zp,
$$ 
where $\mb{A}_F = \mb{A}_F^f \times F_{\infty}$ is the adele ring, $K_f$ is the closure of the image of $\tGa_1(N)$ in $\PGL_2(\mb{A}_F^f)$ (see \ref{K_1}), and $\mc{I}_{\infty}$ is the Iwahori subgroup of $\PGL_2(F_{\infty})$.
For $\phi$ to be special at $\infty$ means that its 
right $\qp[\GL_2(F_{\infty})]$-span is a direct sum of copies of the ``special representation.''  (The latter is the quotient of the locally constant functions $\mb{P}^1(F_{\infty}) \to \qp$ by the constant functions.)

The property of being special at $\infty$ tells us the local behavior at the prime $\infty$ of the $2$-dimensional $\qp$-Galois representation attached to the cusp form.  This is a replacement for the condition of ordinarity at $p$: it is what tells us the $G_{F_{\infty}}$-action on $\mcT$ used in \ref{one_prime}.  

\end{sbpara}

\begin{sbpara}\label{3.7.1}

We explain how the groups $\mcS$ and $\mf{S}$ may be identified.  

The identification passes through the harmonic cocycles on $U(N)$.   These are the functions on the oriented edges of $U(N)$ that change sign if we switch the orientation of an edge and which sum to zero on the edges leading into a vertex (i.e., are harmonic).  The cuspidal harmonic cocycles are those supported on finitely many edges.  The space of $\zp$-valued cuspidal harmonic cocycles may be directly identified with $\mcS$.  It also provides a combinatorial description of $\mf{S}$.  To see this, one starts with the observation that the double coset space on which forms in $\mf{S}$ are defined is none other than the set of oriented edges of $U(N)$.  The property of being special at $\infty$ gives the harmonic condition, and the two notions of cuspidality coincide.  Thus, the spaces $\mf{S}$ and $\mcS$ that appear in the diagram of \ref{key_diag} are canonically identified in the case of $\F_q(t)$.

\end{sbpara}

\subsection{The map $\Upsilon$}

We define the map $\Upsilon \colon Y_{\theta} \to P_{\theta}$ on $\theta$-parts for a fixed primitive character $\theta \colon G \to \qpbt$.

\begin{sbpara}\label{gexact}

We briefly outline the construction of $\Upsilon \colon Y_{\theta} \to P_{\theta}$ that will appear in this section.

As in \ref{Ups_idea}, we analyze the $\mf{h}[G_F]$-action on $\mcT_{\theta}/I_{\theta}\mcT_{\theta}$, showing that it fits in an exact sequence
$$
	0 \to P_{\theta} \to \mcT_{\theta}/I_{\theta}\mcT_{\theta} \to Q_{\theta} \to 0
$$
of $\mf{h}[G_F]$-modules.  Similarly to the setting of $\GL_2$ over $\Q$, the $G_F$-actions on $P_{\theta}$ 
and $Q_{\theta}$ are understood, and $Q_{\theta}$ is free of rank $1$ over $\mf{h}_{\theta}/I_{\theta}$ with a canonical generator.
However, the domain of our map $\Upsilon$ is not a Galois group, so our approach to constructing $\Upsilon$ is different.  We employ compactly supported cohomology, which is dual to Galois cohomology by Poitou-Tate duality.  Instead of directly using the cocycle attached to the exact sequence, we construct $\Upsilon$ in \ref{gup} from a connecting homomorphism $\partial$ on compactly supported \'etale cohomology that appears as the second map in a composition of $\Lambda_{\theta}$-module homomorphisms
$$
	\Upsilon \colon Y_{\theta} \xsim H^2_{\et,c}(\mc{O}[\tfrac{1}{N}],Q_{\theta}(1)) \xrightarrow{\partial} 
	H^3_{\et,c}(\mc{O}[\tfrac{1}{N}],P_{\theta}(1)) 
	\xsim P_{\theta}.
$$ 
The isomorphisms are seen using the $\mf{h}_{\theta}[G_F]$-module structure of $Q_{\theta}$ and the triviality of the $G_F$-action on $P_{\theta}$, respectively.

\end{sbpara}

\begin{sbpara}\label{gxi} 

We define the $L$-function for $\theta$ by
$$
	L(\theta,s) = \prod_{\mf{p} \nmid N} (1 - \theta(\Fr_{\mf{p}})^{-1}\mf{N}(\mf{p})^{-s})^{-1},
$$
where the product is taken over the prime ideals $\mf{p}$ of
$\mc{O}$ not dividing $N$, and $\Fr_{\mf{p}}$ denotes an arithmetic Frobenius at $\mf{p}$.
We then take $\xi_{\theta} \in \Lambda_{\theta}$ to be the nonzero value 
$L(\theta^{-1},-1)$.

\end{sbpara}

\begin{sbpara} \label{h/I-F_q(t)}

We have an isomorphism $\mf{h}_{\theta}/I_{\theta} \xsim \Lambda_{\theta}/(\xi_{\theta})$.  
We indicate one construction of the map.

Consider the Jacobian variety $J$ of $X_1(N)$ and the class $\alpha \in 
J(\C_{\infty})$ of the divisor $(0)-(\infty)$, where $0$ and $\infty$ are cusps on the Drinfeld modular curve.  
Gekeler showed that $\alpha$ has finite order \cite{gekeler}, and it is annihilated by $I$.  The $\mf{h}_{\theta}$-module generated by the $\theta$-part of $\alpha$ is $\Lambda_{\theta}/(\xi_{\theta})$ by a computation of the divisors of Siegel units, providing the desired map.

 \end{sbpara}

\begin{sbpara} \label{Q-F_q(t)}

We define $Q_{\theta} = (\mf{h}_{\theta}/I_{\theta})^{\,\sharp}(1)$, where 
$(\enspace)^{\sharp}$  
indicates a $G_F$-action under which any element that maps to $a \in G$ acts by 
multiplication by $\theta^{-1}(a)$.  Much as in \ref{Q}, 
pairing with the $\theta$-part of $\alpha$ gives rise to a canonical surjection of $\mf{h}_{\theta}[G_F]$-modules $\mc{T}_{\theta}/I\mc{T}_{\theta} \to Q_{\theta}$.

\end{sbpara}

\begin{sbpara} 

The exact sequence
$$
	0 \to P_{\theta} \to \mcT_{\theta}/I_{\theta}\mcT_{\theta} \to Q_{\theta} \to 0
$$
of $\mf{h}[G_F]$-modules is constructed as in \ref{gexact}.  Here, we observe that $Q_{\theta}$
has a nontrivial action of the Frobenius element chosen in \ref{one_prime}, 
so $Q_{\theta}$ is a quotient of 
$\mcT_{\sub}$.  
As before, $\mcT_{\sub}$ is $\zp$-dual to $\mcT_{\quo}$ and 
thereby isomorphic to $\mf{h}$, so we have an isomorphism $\mc{T}_{\sub,\theta}/I_{\theta}\mc{T}_{\sub,\theta}
\xsim Q_{\theta}$ that provides a $G_{F_{\infty}}$-splitting of the exact sequence. 
The known $G_F$-action on $Q_{\theta}$ and the known determinant of the $G_F$-action on
$\mcT_{\theta}$ tell us that $G_F$ acts trivially on $P_{\theta}$.

\end{sbpara}

\begin{sbpara} \label{imc}

The analogue of the Iwasawa main conjecture over $F_N$ is the equality 
$$
	|Y_{\theta}| = [\Lambda_{\theta}:(\xi_{\theta})]
$$ 
of orders.  This equality is a consequence of Grothendieck trace formula, so we do not require the method of Mazur-Wiles to prove it.

\end{sbpara}

\begin{sbpara}\label{gup}
We define our map $\Upsilon$. 

Let $H^i_{\et,c}(\mc{O}[\frac{1}{N}],M)$ denote the $i$th compactly supported \'etale cohomology group of
a compact $\zp\ps{G_F}$-module 
$M$ that is unramified outside $N\infty$.  These groups fit in a long exact sequence
$$
	\cdots \to H^i_{\et,c}(\mc{O}[\tfrac{1}{N}],M) \to H^i_{\et}(\mc{O}[\tfrac{1}{N}],M) \to 
	\bigoplus_{v \mid N\infty} H^i_{\et}(F_v,M) \to H^{i+1}_{\et,c}(\mc{O}[\tfrac{1}{N}],M) \to \cdots.
$$
The exact sequence in \ref{gexact} yields 
a connecting homomorphism 
$$
	H^2_{\et,c}(\mc{O}[\tfrac{1}{N}], Q_{\theta}(1)) \to
	H^3_{\et,c}(\mc{O}[\tfrac{1}{N}], P_{\theta}(1)) = P_{\theta},
$$
the latter identification as $P_{\theta}$ has trivial $G_F$-action, and
we can prove that the canonical map
$$
	H^2_{\et,c}(\mc{O}[\tfrac{1}{N}], Q_{\theta}(1)) \to
	H^2_{\et}(\mc{O}[\tfrac{1}{N}], Q_{\theta}(1))
$$
is an isomorphism.
Hence, the above homomorphism is identified with
$H^2_{\et}(\mc{O}[\tfrac{1}{N}], Q_{\theta}(1)) \to P_{\theta}$.
Our $\Upsilon$ is defined as the following composition:
\begin{equation*}
\begin{split}
	\Upsilon \colon Y_{\theta}  &\xsim
	H^2_{\et}(\mc{O}_N[\tfrac{1}{N}],\Zp(2))_{\theta}\\
	& \xsim
	H^2_{\et}(\mc{O}[\tfrac{1}{N}], \Lambda_{\theta}^{\,\sharp}(2))\\
	& \xsim
	H^2_{\et}(\mc{O}[\tfrac{1}{N}], (\Lambda_{\theta}/\xi_{\theta})^{\,\sharp}(2))\\
	& \xsim H^2_{\et}(\mc{O}[\tfrac{1}{N}], Q_{\theta}(1))
	\to P_{\theta}.
\end{split}
\end{equation*}
The isomorphism in the first line is by \ref{two_obj_F_q(t)}, the isomorphism in the second line is by
Shapiro's lemma, the isomorphism in the third line follows from the fact that $\xi_{\theta}$ kills $Y_{\theta}$ by \ref{imc}, and the isomorphism in the fourth line is by definition of $Q_{\theta}$ in \ref{Q-F_q(t)}.

\end{sbpara}

\subsection{The conjecture: $\varpi$ and $\Upsilon$ are inverse maps}

We state the conjecture and our main result in the case of $\GL_2$ and $\GL_1$ over $\F_q(t)$.

\begin{sbpara}

We state the conjecture.

\begin{conj}
	The maps $\varpi \colon P_{\theta}\to Y_{\theta}$ and $\Upsilon \colon Y_{\theta} \to P_{\theta}$ are
	inverse to each other.
\end{conj}

\end{sbpara}

\begin{sbpara}

We state the theorem.

\begin{thm}\label{th_fun}

We have that $\xi'_{\theta} \Upsilon \circ \varpi=\xi'_{\theta}$, where
$$
	\xi'_{\theta}=\frac{d}{dq^{-s}}L(\theta^{-1}, s)|_{s=-1} \in \Lambda_{\theta}.
$$

\end{thm}

\end{sbpara}

\begin{sbpara}

We can prove the order of $P_{\theta}$ is divisible by the order of $\Lambda_{\theta}/(\xi_{\theta})$ and hence by the order of $Y_{\theta}$.  Thus, in the case that $\xi'_{\theta}$ is a unit in $\Lambda_{\theta}$, our conjecture is implied by the above theorem.

\end{sbpara}

\subsection{The proof that $\xi'\Upsilon \circ \varpi = \xi'$} 

The method of the proof of our Theorem \ref{th_fun} is parallel to the proof in the $\Q$-case.  We give only its bare outline.

\begin{sbpara}\label{newdiag2}  
As in \ref{newdiag}, we consider a refinement of the diagram in \ref{key_diag} in which we
divide the right-hand square of that diagram into two squares:
$$
	\SelectTips{eu}{} 
	\xymatrix{ \mcS_{\theta} \ar[r]^-z \ar[d]^{\bmod I} & H^2_{\et}(Y_1(N), \Z_p(2))_{\theta} 
	\ar@{-->}[r]^-{\HS} \ar[d]^{\infty}
	& H^1_{\et}(\mc{O}_N[\tfrac{1}{N}], \mcT_{\theta}(1)) \ar[r]^-{\text{reg}} \ar[d] & \mf{S}_{\theta} 
	\ar[d]^{\bmod I} \\
	P_{\theta} \ar[r]^{\varpi} & Y_{\theta} \ar[r]^{\xi'_{\theta}}  & Y_{\theta} \ar[r]^{\Upsilon} & P_{\theta}.} 
$$
The commutativity of the leftmost square of the diagram was discussed in Section \ref{zeta_F_q(t)}.  We discuss the maps in the other two squares of the diagram below.

\end{sbpara}

\begin{sbpara}

The map $\HS$ in the diagram arises in a Hochschild-Serre spectral sequence.  An analogue of the discussion
of \ref{HS} applies.

\end{sbpara}

\begin{sbpara}

Let $\kappa$ be the canonical generator of $H^1_{\et}(\F_q, \Z_p)$.
The third vertical arrow in the diagram is the composition
$$
	H^1_{\et}(\mc{O}[\tfrac{1}{N}], \mcT_{\theta}(1))\to H^1_{\et}(\mc{O}[\tfrac{1}{N}], Q_{\theta}(1))
	\xrightarrow{\cup\,\kappa} H^2_{\et}(\mc{O}[\tfrac{1}{N}],  Q_{\theta}(1)) \xsim Y_{\theta},
$$
where  the last isomorphism is given in \ref{gup}.  

\end{sbpara}

\begin{sbpara}\label{gdefp}

The map $\reg$ in the diagram is the $\theta$-part of the $p$-adic regulator map
$$
	H^1_{\et}(F_{\infty}, \mcS(1)) \xrightarrow{\cup \,\kappa} H^2_{\et}(F_{\infty}, \mcS(1)) \xsim
	 \mcS=\mf{S},
$$
where the second map is the invariant map of local class field theory.

\end{sbpara}

\begin{sbpara}

Since $\mcS$ and $\mf{S}$ are canonically identified, both ``mod $I$'' maps in the diagram are just reduction
modulo $I_{\theta}$.

\end{sbpara}

\begin{sbpara}

The proofs of  the commutativity of the other two squares are once again nontrivial, though slightly different, exercises in \'etale and Galois cohomology. 

\end{sbpara}

\begin{sbpara} 

It remains to prove that the composition $\mcS_{\theta}\to \mf{S}_{\theta}\to P_{\theta}$, where the first arrow is the composition of the upper horizontal rows, coincides with $\xi'_{\theta}$ times the reduction modulo $I_{\theta}$ map $\mcS_{\theta}\to P_{\theta}$.  By the computation of Kondo-Yasuda \cite{ky} of the values of a regulator map on the analogues of Beilinson elements, this is reduced to a comparison of their regulator map with the above $p$-adic regulator map. 

\end{sbpara}

\section{What happens for $\GL_d$?} \label{GL_d}

In this section, we discuss three settings for the study of generalizations of the conjectures in Sections 
2 and 3 for $\GL_d$ over a field $F$, for a fixed integer $d \ge 1$.  The fields $F$ and, thereby, the cases we consider here are:
\begin{enumerate}
	\item[(i)] the rational numbers,
	\item[(ii)] an imaginary quadratic field,
	\item[(iii)] a function field in one variable over a finite field.
 \end{enumerate}
We have results only in the cases (i) and (iii) for $d = 2$ discussed above, but we wish to speculate and pose questions in a more general setting.  Rather than formulating precise conjectures, we aim for the more modest goals of pointing in their direction and inspiring the reader to investigate further.

\subsection{The space of modular symbols} \label{symbols}

\begin{sbpara} 

By an infinite place, we mean the unique archimedean place in cases (i) and (ii) and a fixed place $\infty$ in case (iii).  The remaining places are called  finite places.  We have the following objects:
\begin{enumerate}
	\item[$\bullet$] the subring $\mc{O}$ of $F$ of elements that are integral at all finite places,
	\item[$\bullet$] the completion $F_v$ of $F$ at a place $v$,
	\item[$\bullet$] the valuation ring $\mc{O}_v$ of $F_v$ at a nonarchimedean place $v$,
	\item[$\bullet$] the adele ring $\mb{A}_F$ of $F$ and the adele ring $\mb{A}_F^f$ of finite
	places.
	\item[$\bullet$] the subring $\mc{O}_{\mb{A}}^f = \prod_{v \text{ finite}} \mc{O}_v$ of $\mb{A}_F^f$.
\end{enumerate} 
In the discussion below, we will use the notation $(\enspace)^{(d)}$ when defining an object in the
$\GL_d$-setting and then omit the notation in many instances in which $d$ is clear.

\end{sbpara}

\begin{sbpara}

We define a topological space $D_d$ by using the standard maximal compact subgroup 
of $\PGL_d(F_\infty)$: in the respective cases, it is
\begin{enumerate}
	\item[(i)] $\PGL_d(\R)/\PO_d(\R)$, so that $\SL_d(\R)/\SO_d(\R) \xsim D_d$,  
	\item[(ii)] $\PGL_d(\C)/\PU_d$, so that $\SL_d(\C)/\SU_d \xsim D_d$,
	\item[(iii)] the Bruhat-Tits building associated to $\PGL_d(F_\infty)$. 
\end{enumerate}

For example, in case (i) the space $D_2$ is the complex upper half-plane $\mb{H}$.  In case (ii),
the space $D_2$ is the three-dimensional hyperbolic upper-half space $\mb{H}_3$.  Note that in case (iii), the Bruhat-Tits building
has the set $\PGL_d(F_\infty)/\PGL_d(\mc{O}_\infty)$ of homothety classes of $\mc{O}_\infty$-lattices of rank $d$ in $F_{\infty}^d$ as its $0$-simplices.  
\end{sbpara}

\begin{sbpara} \label{K_1}

Let $N$ be a nonzero ideal of $\mc{O}$.
Let $K_1^{(d)}(N)$ be the open compact subgroup of $\GL_d(\mc{O}^f_{\mb{A}})$ given by
$$
	K_1^{(d)}(N) =\Biggl\{g\in \GL_{d}(\mc{O}^f_{\mb{A}}) \mid (g_{d,1}, \ldots, g_{d,d-1}, g_{d,d})
	\equiv (0,\ldots,0,1) \bmod N  \Biggr\}.
$$
Let 
$$
	U^{(d)}(N) =\GL_d(F)\bs (\GL_d(\mb{A}_F^f)/K_1^{(d)}(N) \times D_d).
$$  
The space $U^{(1)}(N)$ is the relative Picard group $\Pic(\mc{O},N)$, viewed as a discrete space.
For $d \ge 2$, the space $U^{(d)}(N)$ is homeomorphic to the disjoint union of  $|\Pic(\mc{O})|$ copies of $\tGa^{(d)}_1(N)\bs D_d$, where 
$$
	\tGa^{(d)}_1(N) = \GL_d(\mc{O}) \cap K_1^{(d)}(N).
$$

\end{sbpara}

\begin{sbpara}\label{epsilon}  Consider cases (i) and (ii).  
Let $\epsilon \in \GL_d(\mc{O})$ be a diagonal matrix with entries $a_1,a_2,\dots,a_d$
such that the product $a_1a_2 \cdots a_d$ generates the roots of unity $\mu_F = \mc{O}^{\times}$ in $F$.
Let 
$$
	\Gamma_1^{(d)}(N) = \tGa_1^{(d)}(N) \cap \SL_d(\mc{O}).
$$
Then  $\tGa_1(N)\bs D_d$ is identified with the quotient of $\Gamma_1(N) \bs D_d$ by the action of the operator 
$$
	\class(g) \mapsto \class(\epsilon g\epsilon^{-1})
$$ 
for $g \in \SL_d(\R)$ in case (i) and for $g\in \SL_d(\C)$ in case (ii).

\end{sbpara}

\begin{sbpara}\label{d=2Q} 
In case (i), the space $U^{(2)}(N)$ is identified with the quotient of $Y_1(N)(\C) = \Gamma_1(N)\bs \mb{H}$ by the action of complex conjugation on $Y_1(N)(\C)$.  In fact, the description in \ref{epsilon} shows that it arises from the quotient of $\mb{H}$ by the action
\begin{multline*}
	\mb{H}\ni x+iy= \class\begin{pmatrix} \sqrt{y}& x \\ 0& 1/\sqrt{y}\end{pmatrix} \\\mapsto  
	\class\left(\begin{pmatrix} -1& 0\\ 0 & 1 \end{pmatrix}
	\begin{pmatrix} \sqrt{y} & x \\0& 1/\sqrt{y}\end{pmatrix} 
	\begin{pmatrix} -1 & 0 \\ 0 & 1\end{pmatrix}\right)
	= \class\begin{pmatrix} \sqrt{y} & -x\\ 0 & 1/\sqrt{y}\end{pmatrix} =-x+iy,
\end{multline*}
which coincides with the action of complex conjugation on $Y_1(N)(\C)$. 
\end{sbpara}

\begin{sbpara} The space $\mcS^{(d)}(N)$ of (cuspidal) modular symbols for $\GL_d$ is defined as
$$
	\mcS^{(d)}(N) = \image(H_{d-1}(U^{(d)}(N), \Z)\to H_{d-1}^{\BM}(U^{(d)}(N), \Z)),
$$
where $H_{*}^{\BM}$ denotes Borel-Moore homology.
Recall that if $\overline{U}(N)$ is a compactification of $U(N)$, 
then $H^{\BM}_i(U(N),\Z)$ is canonically isomorphic to the relative homology group 
$H_i(\overline{U}(N), \overline{U}(N) \setminus U(N), \Z)$.  
The space $\mcS(N)$ may be the homology group $H_{d-1}(\overline{U}(N),\Z)$ for some good choice of compactification.

\end{sbpara}

\begin{sbpara} 

In case (i), we have by \ref{d=2Q} a canonical map
$$
	H_1(X_1(N)(\C),\Z)_+ \to \mcS^{(2)}(N) = H_1(\overline{U}(N),\Z),
$$
where $\overline{U}(N)$ is the quotient of $X_1(N)$ by the action of complex conjugation.
This map is a surjection with $2$-torsion kernel.
 
\end{sbpara}

\begin{sbpara} \label{gen_Hecke}

For a nonzero ideal $\mf{n}$ of $\mc{O}$, let $T(\mf{n})$ denote the Hecke operator on $\mcS^{(d)}(N)$ corresponding to the sum of $K_1^{(d)}(N)$-double cosets of elements of $M_d(\mc{O}_{\mb{A}}) \cap \GL_d(\mb{A}_F^f)$ with determinant generating $\mf{n}\mc{O}_v$ in the $v$-component.  (For $d = 1$, we make the convention that $T(\mf{n}) = 0$ if $\mf{n}$ and $N$ are not coprime.)  These operators satisfy 
$T(\mf{a}\mf{b})=T(\mf{a})T(\mf{b})$ for coprime $\mf{a}$ and $\mf{b}$.  

Let $\mbT^{(d)}(N)$ denote the commutative subring of $\End_{\Z}(\mcS^{(d)}(N))$ generated by the $T(\mf{n})$ with $\mf{n}$ a nonzero ideal of $\mc{O}$.  

\end{sbpara}

\begin{sbpara}

For $d = 1$, the group $\mcS(N)$ of modular symbols is $H_0(U(N),\Z) = \Z[\Pic(\mc{O},N)]$.  The Hecke algebra 
$\mbT(N)$
is the ring $\Z[\Pic(\mc{O},N)]$, with $T(\mf{n})$ for $\mf{n}$ coprime to $N$ equal to the group element for $\mf{n}$.  Under these identifications, $\mbT(N)$ acts by left multiplication on $\mcS(N)$.

\end{sbpara}

\begin{sbpara} The modular symbol $\{0 \to \infty\}$ in Section \ref{modsym} is
generalized to the following element of $ H_{d-1}^{BM}(U(N), \Z)$.
It is the class of the image in the identity component of $U(N)$ of the following standard
subset of $D_d$, with a suitable orientation:
\begin{enumerate}
	\item[(i-ii)]  the set of classes of diagonal matrices in $\GL_d(F_{\infty})$
	with positive real entries,
	\item[(iii)]  the union of all $(d-1)$-simplices with $0$-vertices in
the set of  classes
	in $D_d$ of diagonal matrices in $\GL_d(F_{\infty})$.
\end{enumerate}
The modular symbols $\{\alpha\to \beta\}$ for $\alpha,\beta\in \mb{P}^1(\Q)$ 
are generalized to the classes in $H_{d-1}^{BM}(U(N), \Z)$
of the images in $U(N)$ of the translations by $\GL_d(F)$  of the above
standard subset of $D_d$.

\end{sbpara}

\subsection{Questions for the general case}

We suspect that our results in Sections \ref{GL_2-Q} and \ref{F_q(t)} are special cases of a relationship

\begin{center}

Modular symbols for $\GL_d$ modulo the Eisenstein ideal $\Longleftrightarrow$ Iwasawa theory for $\GL_{d-1}$

\end{center}
that holds for $d \ge 2$.  In this subsection, we describe what we expect to be true.  

\begin{sbpara} We lay out some basic objects, starting with:
\begin{enumerate}
	\item[$\bullet$] a prime number $p \neq \cha F$,
	\item[$\bullet$] a nonzero ideal $N$ of $\mc{O}$ that is coprime to $p$,
	\item[$\bullet$] a commutative pro-$p$ ring $R$ and its total ring of quotients $Q(R)$,
	\item[$\bullet$] a profinite $R$-module $T$ with a continuous
	$R$-linear action of $G_F$ that is unramified at every finite place not dividing $Np$,
\end{enumerate} 
Recalling from 
\ref{gen_Hecke}
the Hecke algebra $\mbT^{(d)}(N)$ and modular symbols $\mcS^{(d)}(N)$, we define
\begin{eqnarray*}
	\mbT_R = \vpr\,(R \otimes \mbT^{(d)}(Np^r)) &\text{and}& \mcS_R = \vpr 
	(R \otimes \mcS^{(d)}(Np^r)).
\end{eqnarray*}
We shall often use the fact that $(p) = \mc{O}$ in case (iii).  For instance, in this case $Np^r = N$, so
we have quite simply that $\mbT_R = R \otimes \mbT^{(d)}(N)$ and $\mcS_R = R \otimes \mcS^{(d)}(N)$.

We also let $\mbT^{(d)}(Np^r)'$ be the subring of $\mb{T}^{(d)}(Np^r)$ generated by the $T(\mf{n})$ with $\mf{n}$ coprime to $(p)$.  Note that $\mb{T}^{(d)}(Np^r)' = \mb{T}^{(d)}(N)$ in case (iii) and $\mb{T}^{(1)}(Np^r)' = \mb{T}^{(1)}(Np^r)$ in
all cases.

\end{sbpara}

\begin{sbpara}\label{RandT}  
We place some conditions on the pair $(R,T)$:
\begin{enumerate}
	\item[(1)] The $Q(R)$-module $V = Q(R) \otimes_R T$ is free of rank $d-1$.

	\item[(2)] For every prime ideal $\mf{p}$ of $\mc{O}$ that does not divide $Np$, 
	the characteristic polynomial 
	$P_{\mf{p}}(u) =  \mr{det}_{Q(R)}(1-\Fr_{\mf{p}}^{-1} u \mid V)$ of an arithmetic Frobenius $\Fr_{\mf{p}}$
	lies in $R[u]$.
\end{enumerate}
For a prime ideal $\mf{p}$ of $\mc{O}$ that does not divide $Np$, 
we define $a(\mf{p}^n)$ for $n \ge 0$ by
$$
	P_{\mf{p}}(u)^{-1} = \sum_{n=0}^\infty a(\mf{p}^n)u^n \in R\ps{u}.
$$
We then suppose:
\begin{enumerate}
	\item[(3)] There exists a ring homomorphism
	$$
		\phi_T \colon \vpr\,(\Z_p \otimes \mbT^{(d-1)}(Np^r)') \to R
	$$ 	
	that sends $T(\mf{p}^k)$ to $a(\mf{p}^k)$ for all prime ideals $\mf{p}$ of $\mc{O}$ not dividing
	$Np$ and all $k \ge 1$.
\end{enumerate}
We extend $a$ to a function on all nonzero ideals $\mf{n}$ of $\mc{O}$ by setting 
$a(\mf{n}) = \phi_T(T(\mf{n}))$ if $\mf{n}$ is coprime to $p$ and $a(\mf{n}) = 0$ otherwise.
In the case $d = 2$, our definition forces $a(\mf{n}) = 0$ for any $\mf{n}$ not coprime to $Np$,
while in general, these values of $a$ may not be uniquely determined by $T$, so $\phi_T$ should
be considered as part of the data.

\end{sbpara}

\begin{sbpara}

We define the Eisenstein ideal $I_T$ of $\mbT_R$ to be the ideal of the $\GL_d$-Hecke algebra $\mbT_R$ generated by the elements
$$
	T(\mf{n})-\sum_{\mf{d}|\mf{n}} a(\mf{d})\mf{N}(\mf{d})
$$ 
for the nonzero ideals $\mf{n}$ of $\mc{O}$.  Note that $I_T$ depends only on $V$ and the choice of $\phi_T$, rather than $T$ itself.  In case (i), the ideal $I_T$ is generated by the coefficients of the formal expression 
$$
	\sum_{n=1}^{\infty} T(n)n^{-s} - \zeta(s) \sum_{\substack{n=1 \\ (n,p) = 1}}^{\infty} a(n)n^{-(s-1)}.
$$

\end{sbpara}

\begin{sbpara}\label{H20} 

For any compact $R\ps{G_F}$-module $M$ that is unramified outside of $S \cup \{\infty\}$ for some finite set $S$ of finite places of $F$ including those dividing $p$, we denote more simply by $H^2_{\et}(\mc{O}[\tfrac{1}{p}], M)$ the $R$-module $H^2_{\et}(\mc{O}[\tfrac{1}{p}], j_*M)$, where $j \colon \Spec(\mc{O})\setminus S
\hookrightarrow \Spec( \mc{O}[\tfrac{1}{p}])$ is the inclusion morphism.  
It is independent of the choice of $S$.  We will also use a similar notation with $\mc{O}$
replaced by its integral closure in a finite extension of $F$.

\end{sbpara}

\begin{sbpara} \label{PY-gen}

Our two objects of study are the $R$-modules:
\begin{enumerate}
	\item[$\bullet$] the geometric object $P= \mcS_R/I_T\mcS_R$ 
	on the $\GL_d$-side,
	\item[$\bullet$] the arithmetic object $Y=H^2_{\et}(\mc{O}[\tfrac{1}{p}], T(d))$ 
	on the $\GL_{d-1}$-side.  
\end{enumerate}
We ask a vague question.

\begin{question}
	Under what conditions does there exist a canonical isomorphism $\varpi \colon P \xsim Y$ of 
	$R$-modules?
\end{question}

We remark that there certainly must be some conditions, as different lattices $T$ in $V$ may have $Y$ that are nonisomorphic.  In what follows, we introduce three settings for further study. 

\end{sbpara}

\begin{sbpara} We fix some notation for abelian extensions of $F$ and their Galois groups.

For $r \ge 0$, let $H_r$ be the ray class field of $F$ of modulus $(p^r)$, and let $O_r$ be the integral closure
of $\mc{O}$ in $H_r$.  Let $\Gamma_r = \Pic(\mc{O}, (p^r))$, which is canonically isomorphic to $\Gal(H_r/F)$ by class field theory.  Let $\Gamma = \varprojlim_r \Gamma_r$.
\begin{enumerate}
	\item[$\bullet$] In case (i), we have that $H_r = \Q(\mu_{p^r})^+$, $O_r = \Z[\zeta_{p^r}]^+$, and 
	$\Gamma = \zp^{\times}/\langle -1 \rangle$.  
	\item[$\bullet$] In case (ii), the field $H_r$ is generated over $F$ by the $j$-invariant $j(E)$ and 
	$x$-coordinates of the $p^r$-torsion points of an elliptic curve $E$ over $F(j(E))$ with CM by $\mc{O}$.  
	There is an exact sequence 
	$$
		0 \to (\zp \otimes \mc{O})^{\times}/\mu_F \to \Gamma \to \Pic(\mc{O}) \to 0.
	$$  
	Note that $\Gamma/\Gamma_{\tor} \cong \zp^2$, where $\Gamma_{\tor}$ is the torsion subgroup of $\Gamma$.
	\item[$\bullet$] In case (iii), we have that $H_r = H_0$, $O_r = O_0$, and
	$\Gamma = \Gamma_r =  \Pic(\mc{O})$.
\end{enumerate}

Let $[a] \in \zp\ps{\Gamma}$ be the group element corresponding to $a \in \Gamma$.  We may also speak of $[\mf{a}]$ for $\mf{a}$ an ideal of $\mc{O}$ coprime to $p$ by taking the sequence of classes of $\mf{a}$ in the groups $\Gamma_r$.
We use $(\enspace)^{\sharp}$ below to denote the (additional) $G_F$-action on a module over a $\zp\ps{\Gamma}$-algebra under which an element that restricts to $a \in \Gamma$ acts by multiplication by $[a]^{-1}$.  

\end{sbpara}

\begin{sbpara}\label{R0T0} We describe setting (A$_d$) for $d \ge 2$.

Let $R_0$ be the valuation ring of a finite extension $K_0$ of $\qp$.  Let $T_0$ be a free $R_0$-module
of rank $d-1$ endowed with a continuous $R$-linear action of $G_F$.  We assume that the $G_F$-action on $T_0$ is unramified at all finite places not dividing $Np$.  We suppose that condition (3) of \ref{RandT} is satisfied for $(R_0,T_0)$, and we use $a_0(\mf{n})$ to denote $a(\mf{n})$ of \ref{RandT} for this pair.

Let $R = R_0\ps{\Gamma}$ and $T=R^{\,\sharp}\otimes_{R_0} T_0$.  
Then the pair $(R,T)$ satisfies conditions (1) and (2) of \ref{RandT}, and we suppose that it satisfies (3).  
It follows directly that $a(\mf{n})=[\mf{n}]^{-1} \otimes a_0(\mf{n})$ for any nonzero ideal $\mf{n}$ of 
$\mc{O}$ that is coprime to $Np$.
By definition of $T$, we also have an $R$-module isomorphism
$$
	Y = H^2_{\et}(\mc{O}[\tfrac{1}{p}], T(d)) \,\cong\, \vpr H_{\et}^2(O_r[\tfrac{1}{p}], T_0(d)).
$$

In that the $G_F$-stable $R_0$-lattice $T_0$ has not been chosen with any special properties inside $V_0 = K_0 \otimes_{R_0} T_0$, we consider an additional condition.  

\begin{enumerate}
	\item[(4)] The $G_F$-representation $k_0 \otimes_{R_0} T_0$ is irreducible over 
	the residue field $k_0$ of $K_0$.
\end{enumerate}
It follows from (4) that the isomorphism class of $T_0$ as an $R_0[G_F]$-module depends only on the $K_0$-representation $V_0$ of $G_F$.  That is, all $G_F$-stable $R_0$-lattices in $V_0$ have the same isomorphism class.  Hence, the isomorphism class of the $R$-module $Y$ depends only on $V_0$.  

Finally, to avoid known exceptions in case (i), we consider a primitivity condition.

\begin{enumerate}
	\item[(5)] The map $\phi_{T_0}$ does not factor through $\vpr (\Z_p \otimes \mbT^{(d-1)}(Mp^r)')$ for
	any ideal $M$ of $\mc{O}$ properly containing $N$.
\end{enumerate}

\end{sbpara}

\begin{sbpara}

We may now ask our question for setting (A$_d$) under conditions (1)--(5).  
\begin{question}  
	Does there exist a canonical isomorphism $\varpi \colon P \xsim Y$ of $R$-modules?
\end{question}
We are also interested in what happens if conditions (4) and (5) are removed.  
For instance, we wonder if (5) might be removed for good choices of $N$, $p$, and $d$, or if (4) might
be removed in the presence of a good, canonical lattice $T_0$.  
In any case, 
we can ask the following question.
\begin{question}
  	If we do not suppose conditions (4) and (5), does there still exist a canonical
	isomorphism $\varpi_{\qp} \colon \qp \otimes_{\zp} P \xsim \qp \otimes_{\zp} Y$?
\end{question}

\end{sbpara}

\begin{sbpara}\label{Cexample} 

The $p$-adic Galois representations $V_0$ attached to the following objects of modulus or level 
$Np^r$ for some $r \ge 0$ all have $(R_0,T_0)$ and $(R,T)$ satisfying (1)--(3):
\begin{enumerate}
	\item[$\bullet$] in case (i) for $d=2$, an even Dirichlet character,
	\item[$\bullet$] in case (i) for $d=3$, a holomorphic cupsidal eigenform,
	\item[$\bullet$] in case (ii) for $d=2$, an algebraic Hecke character on $\mb{A}_F^{\times}$,
	\item[$\bullet$] in case (iii) for $d \ge 2$, a cuspidal eigenform of $\GL_{d-1}$
	that is special at $\infty$.
\end{enumerate}
The examples for $d = 2$ obviously satisfy (4), and in the remaining cases, (4) may be assumed.  
By taking each of the objects to be primitive, we may assume (5).

\end{sbpara}

\begin{sbpara} \label{A_2}

We explain how the setting (A$_2$) for $F = \Q$ and $\F_q(t)$ was studied in Sections \ref{GL_2-Q} and \ref{F_q(t)}.

Let $\theta$ be a primitive character of $\Pic(\mc{O}, Np)$, and impose all the assumptions on $p$, $N$,
and $\theta$ of Sections \ref{GL_2-Q} and \ref{F_q(t)}.
Take $R_0 = \zp[\theta]$, and let $T_0 = \zp[\theta]$ with $G_F$ acting through $\theta^{-1}$.
Let $\Delta = \Pic(\mc{O},(p))$, which we may view as a subgroup of $\Gamma$.  For the objects $P_{\theta}$ and $Y_{\theta}$ of Section \ref{GL_2-Q} in case (i) and of Section \ref{F_q(t)} in case (iii), 
we claim that
\begin{eqnarray*}
	P_{\theta} = \zp \otimes_{\zp[\Delta]} P &\mr{and}& Y_{\theta} = \zp \otimes_{\zp[\Delta]} Y,
\end{eqnarray*}
with $P$ and $Y$ as in \ref{PY-gen}.  This claim is immediate for $\F_q(t)$ as $\Delta$ is trivial, and it is not hard to see for $Y$ in case (i).  However, the claim for $P$ is not evident in case (i), so we prove it. 

\begin{proof}[Proof of the claim]
Note that $R = \zp[\theta]\ps{\Gamma}$ and $T = \zp[\theta]\ps{\Gamma}^{\sharp}$, and note that $\mbT_R = \mbT^{(2)}_R$ of this section is $\mbT[\theta]\ps{\Gamma}$, where $\mb{T}$ is as in \ref{Lambda}.  The claim for $P$ follows if we
can show that the map $T(n) \mapsto T(n)$ on Hecke operators induces an isomorphism
$$
	\mb{T}_{\theta}/I_{\theta} \xsim \zp \otimes_{\zp[\Delta]} (\mbT_R/I_T),
$$
where $I$ is the Eisenstein ideal of \ref{Lambda}.

For a prime $\ell$ not dividing $Np$, the action of $\Fr_{\ell}^{-1}$ on $V$ is multiplication by $\theta(\ell)[\ell]^{-1}$, from which it follows that $a(\ell^k) =  \theta(\ell)^k[\ell]^{-k}$  if $\ell \nmid Np$.
On the other hand, condition (3) forces $a(\ell^k) = 0$ for all $k \ge 1$ for primes $\ell$ dividing $Np$.  The algebra $\mb{T}_R$ contains diamond operators $\langle a \rangle$ for $a \in \Gamma$.  This follows from the identity
$\langle \ell \rangle = \ell^{-1}(T(\ell)^2 - T(\ell^2))$
for $\ell \nmid Np$, which also allows us to compute that
$\langle \ell \rangle \equiv \theta(\ell)[\ell]^{-1} \bmod I_T$.  Thus, $I_T$ is generated by $T(\ell)-1-\ell\langle \ell \rangle$ and 
$\langle \ell \rangle -  \theta(\ell)[\ell]^{-1}$ for primes $\ell \nmid Np$ and $T(\ell)-1$ for primes $\ell \mid Np$.

Noting that the image in $\mb{T}_R/I_T$ of every group element is also the image of an element of $\mb{T}[\theta]$, we now see that the map $T(n) \mapsto T(n)$ induces an isomorphism $(\mb{T}/I)[\theta] \xsim \mbT_R/I_T$ of $\Z[\Delta]$-modules, where $a \in \Delta$ acts by $\theta(a)\langle a \rangle^{-1}$ on
the left and $\theta(a)\langle a \rangle^{-1} \equiv [a] \bmod I_T$ on the right.  The induced map
on $\Delta$-coinvariants is the desired isomorphism.
\end{proof}

\end{sbpara}

\begin{sbpara}

In setting (A$_d$), we have considered Galois cohomology groups of families of $(d-1)$-dimensional 
Galois representations in the variables given by Iwasawa theory.  In case (i) of (A$_d$), for instance, $V$ is a family of Galois representations in the cyclotomic variable.  Of course, there are other families of Galois
representations, such as Hida families, and we would like to consider them.  Therefore, we introduce two additional settings (B$_3$) and (C$_d$) of study.  We do not exclude any representations that are new at $N$ from our families.  
Perhaps we should, but we prefer a simpler presentation.

\end{sbpara}

\begin{sbpara} \label{B}  

We describe setting (B$_3$), in which we work in case (i) for $d = 3$.

Let $\mf{h}$ and $\mcT$ be as in Section \ref{ordinary}, and consider the pair $(\mf{h}^{\circ},\mcT^{\circ}(-1))$, where 
$\circ$ denotes the new-at-$N$ part. 
Condition (1) holds for this pair (see \ref{homology-hecke}). 
As a consequence of Poincar\'e duality, the ordinary \'etale homology group $\mcT$ may be identified with the Tate twist of the ordinary \'etale cohomology group as $\mf{h}[G_F]$-modules.  The characteristic polynomials of $\Fr_{\ell}$ and $T(\ell) \in \mf{h}$ agree on the cohomology $\mcT(-1)$ for any prime $\ell \neq p$.  Thus, condition (2) is satisfied as well, and the map $\phi_{\mc{T}^{\circ}(-1)}$ may be taken to be the identity map on Hecke operators.

Similarly to setting (A$_d$), we consider 
$R= \mf{h}^{\circ}\ps{\Gamma} $ and $T = \mf{h}\ps{\Gamma}^{\,\sharp} \otimes_{\mf{h}} \mcT^{\circ}(-1)$. 
The conditions (1)--(3) are again satisfied for $(R,T)$, and we see that we have $\phi_T$
as in (3) such that $a(n) = T(n)[n]^{-1}$ for $n$ prime to $p$. 
The Eisenstein ideal $I_T$ of $\mbT_R$ is then generated by 
$$
	1 \otimes T(n) - \sum_{\substack{m \mid n\\(m,p) = 1}} mT(m)[m]^{-1}\otimes 1 \in 
	\vpr \mf{h}\ps{\Gamma} \otimes \mbT^{(3)}(Np^r)
$$ 
for all $n\geq 1$.  Note also that we have an $R$-module isomorphism
$$
	Y = H^2_{\et}(\Z[\tfrac{1}{p}], T(3)) \, \cong \, \vpr
	H^2_{\et}(\Z[\zeta_{p^r},\tfrac{1}{p}]^+, \mcT^{\circ}(2)).
$$ 

\end{sbpara}

\begin{sbpara} \label{C}

We describe setting (C$_d$), in which we work in case (iii) for $d \ge 2$.

Let us denote by $Y_1^{(d-1)}(N)$ the Drinfeld modular variety of dimension $d-2$ for $\widetilde{K}_1^{(d-1)}(N)$ over $F$.  We define $T$ by 
$$
	T = \image( H_{c,\et}^{d-2}(Y_1^{(d-1)}(N)_{/\Fbar}, \zp)^{\circ} \to H_{\et}^{d-2}(Y_1^{(d-1)}(N)_{/\Fbar},\zp)^{\circ} ),
$$ 
where $\circ$ denotes the new part (in an appropriate sense).
We then let $R$ be the $\Z_p$-submodule of $\End_{\Z_p}(T)$ generated by the Hecke operators $T(\mf{n})$ for nonzero ideals $\mf{n}$ of $\mc{O}$.  

We imagine but, for $d \ge 4$, are not certain that conditions (1)--(3) hold in this case and that we have $\phi_T$
such that $a(\mf{n}) = T(\mf{n})$ for all $\mf{n}$.  In any case, we may define the Eisenstein ideal $I_T$ of $\mbT_R$ to be generated by 
$$
	1 \otimes T(\mf{n}) - \sum_{\mf{d}|\mf{n}} \mf{N}(\mf{d}) T(\mf{d})\otimes 1
	\in R \otimes \mb{T}^{(d)}(N),
$$ 
for the nonzero ideals $\mf{n}$ of $\mc{O}$.  This Eisenstein ideal is all that we need to consider our question.

\end{sbpara}

\begin{sbpara} 

Our question for (B$_3$) and (C$_d$) is the same as it was for (A$_d$), so we can ask it for all:

\begin{question}  
	Is there a canonical isomorphism $\varpi \colon P \xsim Y$
	of $R$-modules in any of the settings (A$_d$), (B$_3$), or (C$_d$)? 
\end{question}

This question, which has been formulated rather carelessly, is still not fine enough to be a conjecture.  We have more questions than answers: for instance, are the hypotheses that we have made sufficient, and to what extent are they necessary?  What happens for the prime $p = 2$?  We do not wish to exclude it from consideration.  We have made many subtle choices that influence the story in profound yet inapparent ways: e.g., of congruence subgroups, 
Hecke algebras, 
Eisenstein ideals, and \'etale cohomology groups.  Have we made the right choices for a correspondence?  We are glad if the reader is inspired to answer these questions.

\end{sbpara}

\begin{sbpara}\label{hope1}

We end with our hope that it is possible to explicitly define the maps $\varpi$ that are the desired isomorphisms in the settings (A$_d$), (B$_3$), and (C$_d$).

The groups $\mcS(N)$ often have explicit presentations very similar to those of Section \ref{modsym}.  These are found in the work of Cremona \cite{cremona}, Ash \cite{ash}, Kondo-Yasuda \cite{ky}, and others.  So, explicit definitions of $\varpi$ and affirmative answers to our questions would give explicit presentations of the arithmetic object $Y$. 

The map $\varpi$ should take a modular symbol to a cup product of $d$ special units.  As explained above, this has been done in cases (i) and (iii) for $d = 2$.  Beyond these, the settings in which we hope to do this are:
\begin{enumerate}
	\item[$\bullet$] (A$_2$) in case (ii), using cup products of two elliptic units,
	\item[$\bullet$] (B$_3$) using cup products of three Siegel units, 
	\item[$\bullet$] (C$_d$) using cup products of $d$ of the Siegel units in \cite{ky}.
\end{enumerate}	
Goncharov has made closely related investigations into the first two of these settings \cite{goncharov}.

\end{sbpara}

\vspace{1ex}

\begin{ack}

{The work of the first two authors (resp., third author) was supported in part by the National Science Foundation under Grant Nos.~DMS-1001729 (resp., DMS-0901526).}

\end{ack}

\renewcommand{\baselinestretch}{1}

\end{document}